\magnification=\magstep1
\nopagenumbers \headline={\tenrm\hfil --\folio--\hfil}
\baselineskip=18pt  \lineskip=3pt minus 2pt   \lineskiplimit=1pt
\hsize=15true cm  \vsize=24true cm 
\def\tomb{\phantom{.}\hfill\vrule height.4true cm width.3true cm \par\smallskip\noindent}
\def\ker{{\rm Ker\,}}  \def\chop{\hfill\break}  
\def\f #1,#2.{\mathsurround=0pt \hbox{${#1\over #2}$}\mathsurround=5pt}
\def\s #1.{_{\smash{\lower2pt\hbox{\mathsurround=0pt $\scriptstyle #1$}}\mathsurround=5pt}}

\def\mapdownl #1;{\vcenter{\hbox{$\scriptstyle#1$}}\Big\downarrow}
\def\mapdownr #1;{\Big\downarrow\rlap{$\vcenter{\hbox{$\scriptstyle#1$}}$}}
\tolerance=1600 \mathsurround=5pt  
\def\maprightu #1;{\smash{\mathop{\longrightarrow}\limits^{#1}}}
\def\maprightd #1;{\smash{\mathop{\longrightarrow}\limits_{#1}}}
\def\maprightt #1,#2.{\mathrel{\smash{\mathop{\longrightarrow}\limits_{#1}^{#2}}}}

\def\brc #1,#2.{\left\langle #1\,|\,#2\right\rangle}
\def\rn#1{{\romannumeral#1}}  \def\cl #1.{{\cal #1}}
\def\convl #1,#2.{\mathrel{\mathop{\longrightarrow}\limits^{#1}_{#2}}}
\def\convr #1,#2.{\mathrel{\mathop{\longleftarrow}\limits^{#1}_{#2}}}
\def\set #1,#2.{\left\{\,#1\;\bigm|\;#2\,\right\}}
\def\theorem#1@#2@#3\par{\smallskip\parindent=.6true in \itemitem{\bf #1}
{\sl #2}\parindent=20pt\smallskip\itemitem{\it Proof:\/}#3\tomb}
\def\thrm#1"#2"#3\par{\smallskip\parindent=.6true in \itemitem{\bf #1}
{\sl #2}\parindent=20pt\smallskip\itemitem{\it Proof:\/}#3\tomb}
\def\teorem#1@#2@ {\smallskip\parindent=.6true in \itemitem{\bf Theorem #1}
{\sl #2}\hfill \parindent=20pt\smallskip\noindent} \def\ab{\allowbreak}
\def\fteorem#1@#2@ {\smallskip\parindent=.6true in \itemitem{\bf #1}{\sl #2}
\hfill\parindent=20pt\smallskip\noindent}

\def\shift #1;{\mathord{\phantom{#1}}}  
\def\ref #1.{\mathsurround=0pt${}^{#1}\phantom{|}$\mathsurround=5pt}
\def\cross #1.{\mathrel{\jackup 3,\mathop\times\limits_{#1}.}}
\def\jackup#1,#2.{\raise#1pt\hbox{\mathsurround=0pt $#2$\mathsurround=5pt}}

\def\bbrc #1,#2,#3.{\langle #1 |\,#2\,|#3\rangle}
   
\def\cst#1,{{C^*(#1-\un)}}
\def\un{\hbox{\mathsurround=0pt${\rm 1}\!\!{\rm 1}$\mathsurround=5pt}}

\def\alg#1.{{C^*(#1)}}
\def\slim{\mathop{{\rm s\!\!-\!\!lim}}}
\def\ccr #1,#2.{{\overline{\Delta(#1,\,#2)}}}
\input amssym.def
\input amssym.tex
\def\un{\Bbb I}

\def\C{\Bbb C}
\def\b#1'{{\bf #1}}
\def\rest{\mathord{\restriction}}

\def\rep{{\rm Rep}\,}

\def\sp{{\rm Span}\,}

\def\wt{\widetilde}
\def\wh{\widehat}
\def\ol#1.{\overline{#1}}
\def\rsg{{\rm Rep}_\sigma G}
\def\rbg{{\rm Rep}^B_\sigma G}
\def\rsd{{\rm Rep}_\sigma G_d}
\def\csd{{\rm C}^*_\sigma(G_d)}
\def\csg{{\rm C}^*_\sigma(G)}
\def\sga{$\sigma\hbox{--group}$algebra }

\pageno=1  \noindent
\centerline{\bf Group Algebras for Groups which are not Locally Compact.}
\vglue .5in
\centerline{Hendrik Grundling } 
\centerline{ School of Mathematics, University of
New South Wales,}
\centerline{ Sydney, NSW 2052, Australia.}
\centerline{ email: hendrik@maths.unsw.edu.au} 
\vglue .5in
\item{{\bf Abstract}}{\sl 
We generalise the definition of a group algebra so that
it makes sense for  general topological groups,
not necessarily locally compact. In
particular, we require that the representation theory
of the group algebra is isomorphic (in the sense of Gelfand--Raikov)
to the continuous representation theory of the group,
or to some other important subset of representations.
We prove that a group algebra if it exists, is always unique
up to isomorphism. From examples, group algebras do not always
exist for non--locally compact groups, but they do exist for some. 
We define a convolution on the dual of the Fourier--Stieltjes 
algebra making it into a C*-algebra, we prove that a group algebra
 if it exists, can always be
embedded in this convolution algebra, 
and we find sharp conditions for a subalgebra to be a group
algebra. When the group is locally compact, we obtain a new characterisation
of its group algebra which does not involve the Haar measure, nor behaviour of
measures on compact sets.}
\chop
\chop {\bf AMS classification: 43A10, 43A20, 43A35, 46L05}

\beginsection  Introduction.

The Gelfand--Raikov theorem [GR]
proved that the continuous (unitary) representation
theory of any locally compact group is isomorphic in a natural sense
to the (nondegenerate Hilbert space) representation theory of a C*--algebra. 
The proof is constructive,
in that it explicitly constructs the group algebra as the enveloping
C*-algebra of the convolution algebra $L^1(G),$
and faithfully embeds the group as unitaries in
the multiplier algebra of the group algebra. Subsequently group algebras
for locally compact groups have been generalised in many directions
(e.g. twisted group algebras, groupoid algebras, some semigroup algebras
and cross--products of a C*-algebra by a group action), and has been a central component of
harmonic analysis. The question naturally arises as to whether the
Gelfand--Raikov theorem can be extended to topological groups which are 
NOT locally compact, and this problem will be at the focus of our investigation here.
This question has attained some urgency, due to the fact that
such groups regularly and naturally arise in physics and mathematics
(e.g. gauge groups, diffeomorphism groups, symplectomorphism groups,
Banach Lie groups, inductive limit groups, Fr\'echet-Lie groups etc.)
and a substantial body of work is currently developing on the continuous representation
theory of these groups. Group algebras are very useful to have, e.g. they allow one
to use the topological analysis of the spectrum of C*--algebras
to analyze group representations,
use direct integral decomposition theory for representations
of (separable) C*--algebras to decompose continuous group representations into
other continuous group representations, and
use Rieffel induction to induce group representations from subgroups
to larger groups.

In its full generality, the question has a negative answer, i.e. it is not
true for all topological groups that their continuous representation theory 
is isomorphic (in the sense of Gelfand--Raikov) to the representation theory
of a C*--algebra. For instance, there are Abelian groups with NO nontrivial
continuous
unitary representations, cf. Banaszczyk~[Ban], and there are Abelian groups
with continuous representations, but no irreducible ones cf. Example~5.2 in Pestov~[Pes].
However, we can still ask the question of whether some useful subset of representations
of a topological group is isomorphic to the representation theory of a C*--algebra.
Of course one would also like to characterize those topological groups for which the
Gelfand--Raikov theorem holds, i.e. for which their continuous unitary representation
theory is isomorphic to the representation theory of a C*--algebra.

To make the discussion more precise, let us fix notation.
For a C*-algebra $\cl A.$ denote its set of nondegenerate Hilbert space representations
by $\rep\cl A..$
For a topological group $G$ with a fixed nondiscrete Hausdorff topology, let
$\sigma:G\times G\to{\Bbb T}$ be a 2--cocycle $\sigma\in
Z^2(G,\,{\Bbb T})$ which is jointly continuous and normalised, i.e.
$1=\sigma(x,\,e)=\sigma(e,\, x)=\sigma(x,\,x^{-1})$ and
$\ol\sigma(x,\,y).=\sigma(y^{-1},\,x^{-1})$ for all $x,\,y\in G.$
Let $\rsg$ denote the set of strong operator continuous unitary 
$\sigma\hbox{--representations}$ on Hilbert space, and in the case
$\sigma=1,$ we simplify the notation to $\rep G.$
Let $G_d$ denote $G$ with its discrete topology. 
Clearly $\rsd\big(\supseteq\rsg\big)$ is the set of {\it all} unitary
$\sigma\hbox{--representations}$ of $G,$ not necessarily continuous,
and this is isomorphic to to $\rep\csd$ where $\csd$ is the group
algebra of $G_d.$

Structurally, continuous group representation theory is quite similar to
C*-algebra representation theory, but there are also differences.
Note first that both of $\rep G$ and $\rep\cl A.$ are:
\item{$\bullet$} closed w.r.t. composition with (continuous) homomorphisms,
\item{$\bullet$} closed with respect to direct sums of representations,
\item{$\bullet$} closed with respect to unitary conjugation,
\item{$\bullet$} closed with respect to subrepresentations, i.e. if for a representation
$\pi$ on $\cl H._{\pi}$ there is a closed invariant subspace $\cl K.\subset\cl H._{\pi},$
then $\pi(G)$ (resp. $\pi(\cl A.))$ restricted to $\cl K.$
produces a representation in $\rep G$ (resp. $\rep\cl A.).$

\noindent As for the differences, note that $\cl A.$ is separated by its
irreducible representations, but the continuous irreducible representations
 need not separate $G.$ Secondly,
under tensor products $\rep G$ is closed (or more generally,
$\rsg\otimes{\rm Rep}_\rho G\subseteq{\rm Rep}_{\sigma\rho} G)$
but the tensor product of two representations of $\cl A.$ is a representation
of $\cl A.\otimes\cl A.,$ not of $\cl A..$
Our concept of isomorphism between $\rsg$ and $\rep\cl A.$ is:
 \vfill\eject
\item{\bf Def.} 
Let $G$ be a topological group, $\sigma$ as above,
and let $\cl R.\subset\rsd$ be a given subset of unitary $\sigma\hbox{--representations}$
of $G.$ Then a
   $\sigma\hbox{--\bf group}$ {\bf algebra}
for the pair $(G,\,\cl R.)$
is a C*--algebra $\cl L.$ and a $\sigma\hbox{--homomorphism}$
$\varphi:G\to UM(\cl L.)$
such that the unique extension map $\theta:{\rm Rep}\cl L.\to\rsd$
is injective, and with image $\theta\big({\rm Rep}\,\cl L.\big)=\cl R..$\chop
In this case we say that $\cl R.$ is {\bf isomorphic} to $\rep\cl L..$
\itemitem{\bf Remark.}(1) The map $\varphi$ maps to the unitaries in the multiplier algebra of
$\cl L.$ and satisfies $\varphi(g)\varphi(h)=\sigma(g,h)\varphi(gh).$
Moreover, any nondegenerate representation of $\cl L.$ has a unique extension to its
multiplier algebra $M(\cl L.),$ and this defines the map $\theta:{\rm Rep}\cl L.\to\rsd$
by $\theta(\pi)(g):=\slim\limits_{\alpha\to\infty}\pi\big(\varphi(g)E_\alpha\big)$
where $\{E_\alpha\}\subset\cl L.$ is any approximate identity of $\cl L..$
Note that we may have that $\varphi$ is not injective, e.g. in the case
when $\cl R.$ does not separate $G.$ 
\item{(2)}
When  $G$ is locally compact, the usual $\sigma\hbox{--group algebra}$ will be for
the case that $\cl R.=\rsg,$ and then $\cl L.=\csg$ satisfies the
conditions of the definition, where $\varphi$ is injective. Below we will prove uniqueness.
This generalisation of group algebras seem useful even for locally compact groups,
because it allows the analysis of representation sets other than
$\rsg.$
\item{(3)} 
For a small class of non-locally compact groups, $\sigma\hbox{--group algebras}$
were constructed for $\cl R.=\rsg$ in [Gr1].
The existence question for group algebras was studied in a more general
context in~[Gr2].
\item{(4)} 
Note that the map $\theta$ preserves direct sums, unitary conjugation, subrepresentations,
and (as we will see) irreducibility, so that this notion of isomorphism between
$\cl R.$ and $\rep\cl L.$ involves strong structural correspondences, and restricts
the class of sets $\cl R.$ for which group algebras exist.
However, this isomorphism is obviously not an equivalence relation, since it
relates objects in two distinct sets.
In the case that $\theta:{\rm Rep}\cl L.\to\cl R.$ is surjective but not injective,
it is natural to say that $\rep\cl L.$ is {\it homomorphic} to  $\cl R.,$ since $\theta$
still transfers some structure to $\cl R.$ (though e.g. irreducibility
of representations is lost).
We will not examine this concept here.

\noindent In the rest of this paper, we will develop the concept of group algebras
for groups which are not locally compact. Below, in Sect.~1 we consider general
theory and prove that group algebras, if they exist, are unique up to isomorphism.
In the subsequent sections, we will be concerned with the existence question
when $\cl R.$ is a subset of the continuous representations $\rsg,$ and in particular
we will
analyze an important convolution algebra, and develop 
conditions to ensure
the existence of a group algebra.

\beginsection 1. General structures for group algebras.

Let $\cl L.$ be a C*-algebra, and recall that the {\bf strict topology}
of its multiplier algebra $M(\cl L.)$is given by the family of seminorms on $M(\cl L.):$
$$B\to\|BA\|+\|AB\|,\quad A\in\cl L.,\;B\in M(\cl L.)\;.$$
Then $\cl L.$ is strictly dense in $M(\cl L.),$ cf. Prop.~3.5 and 3.6 in~[Bus].
Below $\sp X$ will denote the space of finite linear combinations of $X.$
\thrm Proposition 1.1."Let $G,\;\sigma$ as above, let
 $\cl R.\subset\rsd$ be given,  and let
$\cl L.$ be a $\sigma\hbox{--group}$ algebra for $(G,\,\cl R.).$ Then\chop
{\bf (1)} $\sp\varphi(G)$ is a strictly dense *--algebra in $M(\cl L.),$\chop
{\bf (2)} Each $\pi\in\rep\cl L.$ is strict--strong operator continuous, and
   $\theta(\pi)$ is the strict extension of $\pi$ to $\varphi(G).$\chop
{\bf (3)} Each $\pi\in\cl R.$ has a unique extension to  $\sp\varphi(G)$ 
  which is strict--strong operator continuous, and conversely
  $\cl R.$ is exactly the restrictions to $\varphi(G)$ of the
  set of strict--strong operator continuous representations of
   $\sp\varphi(G).$\chop
{\bf (4)} The inverse map of the bijection $\theta,$ is the map
 $\theta^{-1}:\cl R.\to\rep\cl L.$ obtained by
 $\theta^{-1}(\pi)(A):=\slim\limits_\alpha\wt\pi(B_\alpha)$
 where $\wt\pi$ is the unique strictly--strong operator continuous extension
 in (3), and $\{B_\alpha\}\subset\sp\varphi(G)$ is a net strictly converging
 to $A\in\cl L..$\chop
{\bf (5)} If $G$ is Abelian and $\sigma=1,$ then $\cl L.$ is
commutative."
(1) That  $\sp\varphi(G)$ is a *--algebra is obvious from the fact that
$\varphi$ is a $\sigma\hbox{--homomorphism}.$ Let
$\cl Q.$ be the strict closure of  $\sp\varphi(G).$ 
This is a *-algebra, so since $\varphi(G)$ separates $\cl R.=
\theta(\rep\cl L.),$ it follows that $\cl Q.$ separates 
$\rep\cl L..$ Thus by Prop.~2.2 in [Wor], we have that
$\cl Q.=M(\cl L.).$\chop
(2) Let $\pi\in\rep\cl L.,$ which is a *--homomorphism
$\pi:\cl L.\to\pi(\cl L.)=:\cl C.\subset\cl B.(\cl H._\pi),$
and by Prop.~3.8 and 3.9 in [Bus], this extends uniquely to a *--homomorphism
$\pi:M(\cl L.)\to M(\cl C.)\subseteq\cl B.(\cl H._\pi)$
which is strict--strict continuous (using nondegeneracy of $\pi).$
Since on $M(\cl C.)\subseteq\cl B.(\cl H._\pi)$ the strong operator topology is
coarser than the strict topology, it follows that 
$\pi:M(\cl L.)\to\cl B.(\cl H._\pi)$ is strict--strong operator continuous.
If $\{E_\alpha\}\subset\cl L.$ is an approximate identity of $\cl L.,$
then for each $B\in M(\cl L.)$ the net $\{BE_\alpha\}$ strictly converges
to $B,$ hence $\pi(BE_\alpha)$ converges in strong operator topology
to $\pi(B),$ and by definition this is $\theta(\pi)(g)$ when $B=\varphi(g).$\chop
(3) By the bijection $\theta:\rep\cl L.\to\cl R.,$ for each $\pi\in\cl R.$
there is a $\rho\in\rep\cl L.$ such that its strict extension to $M(\cl L.)$
produces $\pi\in\cl R.$ by (2). 
Hence each $\pi\in\cl R.$ has a strictly continuous extension
$\wt\pi=\rho\rest\sp\varphi(G)$ to $\sp\varphi(G).$
If $\wh\pi$ is another strictly continuous extension of $\pi$
to $\sp\varphi(G),$ then since $\sp\varphi(G)$ is strictly dense in
$M(\cl L.),$ it extends uniquely to $\cl L.,$ so by definition we get
$\theta(\wt\pi)=\pi=\theta(\wh\pi).$
Since $\theta$ is injective, we have that 
$\wt\pi\rest\cl L.=\wh\pi\rest\cl L.$ and as $\cl L.$ is strictly dense
we have that $\wh\pi=\wt\pi\,.$
Conversely, if $\pi$ is a (bounded) *-representation of $\sp\varphi(G)$
which is strictly continuous, then it extends uniquely to $M(\cl L.),$
in which case $\theta\big(\pi\rest\cl L.\big)=\pi\rest\varphi(G)
\in\cl R..$\chop
(4) This is clear from the previous parts.\chop
(5) The strict topology on $M(\cl L.)\subset\cl L.''$ is finer than the
weak operator topology of $\cl L.''$  on $M(\cl L.).$
Thus the strict closure of $\sp\varphi(G)$ (i.e. $M(\cl L.))$ is contained 
in its weak operator closure, and this is the double commutant
$\sp\varphi(G)''$ since $\varphi(G)$ contains the identity.
Now if $G$ is Abelian and $\sigma=1$ we have that $\sp\varphi(G)$
is commutative, hence $\sp\varphi(G)''\supset\cl L.$ is commutative.

Since $\sp\varphi(G)$ and $\cl L.$ are both strictly dense in $M(\cl L.),$
and the strictly continuous representations on $M(\cl L.)$ are the 
extensions of representations in $\rep\cl L.,$ it follows that all
the properties of these representations are determined by their restrictions
to either $\varphi(G)$ or $\cl L..$ Below we will use this to take sufficient
analytical structure from $\rep\cl L.$ to $\cl R.,$ so that we can characterize
$\cl L.$ directly on $\cl R.$ in our uniqueness theorem.

Let $\cl L.$ be a C*-algebra, then we topologize $\rep\cl L.$ according to
Takesaki and Bichteler~[Bic,~Tak]. Let $\pi\s{\cl L.}.$ be the universal representation
of $\cl L.$ on the universal Hilbert space $\cl H..$
Then we can identify $\rep\cl L.$ with the set of all (including degenerate)
representations of $\cl L.$ on $\cl H.,$ and we denote this by
$\rep(\cl L.,\,\cl H.).$ Then we equip $\rep(\cl L.,\,\cl H.)$
with the pointwise strong operator topology, i.e. a net
$\{\pi_\nu\}\subset\rep(\cl L.,\,\cl H.)$
converges to $\pi\in\rep(\cl L.,\,\cl H.)$ iff
$\pi_\nu(A)\to\pi(A)$ in the strong operator topology of $\cl B.(\cl H.)$
for all $A\in\cl L..$

From Prop.~1.1 we obtain:
\fteorem Corollary 1.2.@
Let 
$\cl L.$ be a $\sigma\hbox{--group}$ algebra for $(G,\,\cl R.).$ Then
$\pi\in\rep\cl L.$ is cyclic (resp. irreducible) iff 
$\theta(\pi)\in\cl R.$ is cyclic (resp. irreducible). @
From this fact we get that 
$$\theta(\pi\s{\cl L.}.)=\oplus\set\pi\in\cl R.,\pi\;\hbox{has a cyclic vector $\Omega$
with $\|\Omega\|=1$}.$$
since $\pi\s{\cl L.}.$ is the direct sum of GNS--representations of the states
of $\cl L..$
This characterises $\theta(\pi\s{\cl L.}.)$ directly on $G,$ so we can consider
$\cl R.\subset\rep_\sigma(G_d,\,\cl H.)\equiv\sigma\hbox{--representations}$
 of $G$ on $\cl H.$ 
which are unitary on their essential subspaces,
and hence we topologise $\cl R.$ also with the
pointwise (on $G)$ strong operator topology. From Prop.~1.1 we also get:
\fteorem Corollary 1.3.@ Let $\cl L.$ be a $\sigma\hbox{--group}$ algebra 
for $(G,\,\cl R.),$ and $\rep(\cl L.,\,\cl H.)$ as above. Denote the
essential subspace of $\pi\in\rep(\cl L.,\,\cl H.)$ by $\cl H._\pi\,,$
with essential projection $P_\pi:\cl H.\to\cl H._\pi.$
\chop
{\bf (1)} If $\pi,\,\pi'\in\rep(\cl L.,\,\cl H.)$ satisfy $\cl H._\pi\perp
\cl H.\s{\pi'}.,$ then $\theta(\pi\oplus\pi')=\theta(\pi)\oplus\theta(\pi')\,.$
Conversely, if $\pi,\;\pi'\in\cl R.\subset\rep_\sigma(G_d,\,\cl H.)$
with $\cl H._\pi\perp\cl H.\s{\pi'}.,$ then 
$\theta^{-1}(\pi\oplus\pi')=\theta^{-1}(\pi)\oplus\theta^{-1}(\pi')\,.$\chop
{\bf (2)} The essential projections of $\pi\in\rep(\cl L.,\,\cl H.)$
and $\theta(\pi)\in\cl R.\subset\rep_\sigma(G_d,\,\cl H.)$
are the same\chop
{\bf (3)} for $\pi\in\rep(\cl L.,\,\cl H.)$ let $U\in\cl B.(\cl H.)$ be any 
partial isometry with initial projection $U^*U\geq P_\pi,$ and define
$\pi^U(A):=U\pi(A)U^*,$ $A\in\cl L..$ Then
$\theta(\pi^U)=U\theta(\pi)U^*=:\theta(\pi)^U\,,$ and conversely
$\theta^{-1}(\pi^U)=\theta^{-1}(\pi)^U$ for all $\pi\in\cl R..$  @
Thus by these two corollaries, $\theta$ preserves much of the structure of
$\rep(\cl L.,\,\cl H.).$ In fact, it also preserves the topology:
\thrm Proposition 1.4."
 Let $\cl L.$ be a $\sigma\hbox{--group}$ algebra
for $(G,\,\cl R.),$ then\chop
 $\theta:\rep(\cl L.,\cl H.)\to\cl R.$
is a homeomorphism."
 From Prop.~1.1, it suffices to show that for the
strict extensions of$\rep(\cl L.,\,\cl H.)$ to $M(\cl L.),$
for a net $\{\pi_\nu\}$ we have the convergence
$\pi_\nu(A)\to\pi(A)$ for all $A\in\cl L.$ in strong operator topology
iff $\pi_\nu(B)\to\pi(B)$ for all $B\in\varphi(G)$ in strong operator topology.
\chop
Assume that $\pi$ and $\{\pi_\nu\}$ are  strict--strong operator continuous
representations in $\rep(M(\cl L.),\cl H.)$ such that
$\pi_\nu(A)\to\pi(A)$ for all $A\in\cl L.$ in strong operator topology.
For any $B\in\varphi(G),$ let
$\{A_\alpha\}\subset\cl L.$ be a net strictly converging to
$B.$ Then for all $\psi\in\cl H.$ we have:
$$\eqalignno{\big\|(\pi_\nu(B)-\pi(B))\psi\big\|&\leq
\big\|\pi_\nu(B-A_\alpha)\psi\big\| +
\big\|(\pi_\nu(A_\alpha)-\pi(A_\alpha))\psi\big\| \cr
&\qquad\qquad +\big\|\pi(A_\alpha-B)\psi\big\|\;.&-(1)\cr}$$
Since $A_\alpha\to B$ strictly, we have for the extension of
the universal representation  $\pi\s{\cl L.}.$ to $M(\cl L.)$ that
${\big\|\pi\s{\cl L.}.(B-A_\alpha)\psi\big\|}\to 0$ for all
$\psi\in\cl H..$ Since $\pi$ and $\pi_\nu$ are subrepresentations
of $\pi\s{\cl L.}.,$ we also have that
$$\big\|\pi_\nu(B-A_\alpha)\psi\big\|\leq
\big\|\pi\s{\cl L.}.(B-A_\alpha)\psi\big\|
\geq \big\|\pi(B-A_\alpha)\psi\big\|\,.$$
Thus for each $\varepsilon>0$ there is an $\alpha_1$ such that
for all $\nu$
$$\big\|\pi_\nu(B-A_\alpha)\psi\big\|+
\big\|\pi(B-A_\alpha)\psi\big\|\leq\varepsilon\quad\forall\;
\alpha>\alpha_1\,.$$
Thus from $(1)$ we get for all $\alpha>\alpha_1$ that
$$\eqalignno{
\lim_{\nu\to\infty}\big\|(\pi_\nu(B)-\pi(B))\psi\big\|&\leq
\varepsilon+\lim_{\nu\to\infty}\big\|(\pi_\nu(A_\alpha)-\pi(A_\alpha))\psi\big\|\cr
&=\varepsilon\;.\cr}$$
So, since $\varepsilon>0$ is arbitrary, we have for all $\psi\in\cl H.$ that\chop
$\lim\limits_{\nu\to\infty}\big\|(\pi_\nu(B)-\pi(B))\psi\big\|=0\;.$\chop
Conversely, assume that for  $\pi$ and $\{\pi_\nu\}$ strict--strong operator continuous
representations in $\rep(M(\cl L.),\cl H.),$ that
$\pi_\nu(B)\to\pi(B)$ for all $B\in\varphi(G)$ in strong operator topology.
By triangle inequalities, we then get that
$\pi_\nu(B)\to\pi(B)$ for all $B\in\sp\varphi(G)$ in strong operator topology.
For any $A\in\cl L.,$ let $\{B_\alpha\}\subset\sp\varphi(G)$
be a strictly convergent net to $A\in\cl L..$ As above, we get that for
any $\varepsilon>0$ and $\psi\in\cl H.$ there is an $\alpha_1$ such that
for all $\nu$
$$\leqalignno{\big\|\pi_\nu(A-B_\alpha)\psi\big\|&+
\big\|\pi(A-B_\alpha)\psi\big\|\leq\varepsilon\quad\forall\;
\alpha>\alpha_1\,.\cr
{}\qquad\quad\big\|(\pi_\nu(A)-\pi(A))\psi\big\|&\leq
\big\|\pi_\nu(A-B_\alpha)\psi\big\| +
\big\|(\pi_\nu(B_\alpha)-\pi(B_\alpha))\psi\big\|&\hbox{Thus:} \cr
&\qquad\qquad +\big\|\pi(B_\alpha-A)\psi\big\|\cr
&\leq\varepsilon+\big\|(\pi_\nu(B_\alpha)-\pi(B_\alpha))\psi\big\|\cr}$$
for all $\alpha>\alpha_1\,.$ Take the limit $\nu\to\infty$ on both sides,
and use the fact that $\varepsilon>0$ is arbitrary to find that
$\lim\limits_{\nu\to\infty}\big\|(\pi_\nu(A)-\pi(A))\psi\big\|=0\;.$

Following 
Takesaki and Bichteler~[Bic,~Tak], we define:
\item{\bf Def.} An {\bf admissible operator field}
on $\rep(\cl L.,\cl H.)$ (resp. $\cl R.)$ is a map\chop
$T:\rep(\cl L.,\cl H.)\to\cl B.(\cl H.)$
(resp. $T:\cl R.\to\cl B.(\cl H.))$ such that:
\item{(\rn1)}
$\|T\|=\sup\set\|T(\pi)\|,{\pi\in\rep(\cl L.,\cl H.)}.<\infty$\chop
(resp. $\|T\|=\sup\set\|T(\pi)\|,{\pi\in\cl R.}.<\infty),$
\item{(\rn2)}
$T(\pi)=P_\pi T(\pi)=T(\pi)P_\pi$ for all $\pi\in\rep(\cl L.,\cl H.)$
(resp. $\pi\in\cl R.)$ where $P_\pi$ denotes the essential projection of $\pi,$
\item{(\rn3)}
$T(\pi\oplus\pi')=T(\pi)\oplus T(\pi')$ whenever $\cl H._\pi\perp
\cl H._{\pi'}$ in $\cl H.,$
\item{(\rn4)}
$T(\pi^U)=UT(\pi)U^*$ for all $\pi\in\rep(\cl L.,\cl H.)$
(resp. $\pi\in\cl R.)$ where
$U\in\cl B.(\cl H.)$ is a partial isometry with
$U^*U>P_\pi\,.$

\noindent
The set of admissible operator fields form a C*--algebra under pointwise 
operations and the sup--norm, and we denote the two resultant C*-algebras by
${\cl A.(\cl L.,\cl H.)}$ and ${\cl A.(\cl R.,\cl H.)}$
respectively. In particular ${\cl A.(\cl L.,\cl H.)}$ contains 
the C*-algebra
$$\wt{\cl L.}:=\set {T_A:\rep(\cl L.,\cl H.)\to\cl B.(\cl H.),\;
A\in\cl L.}, {T_A(\pi):=\pi(A)\quad\forall\;\pi\in\rep(\cl L.,\cl H.)}.$$
which is obviously isomorphic to $\cl L..$
Then we have the Takesaki--Bichteler duality theorem:
$$\wt{\cl L.}\cong\set {T\in\cl A.(\cl L.,\cl H.)},
T\quad\hbox{is strong--operator continuous}.$$
where $\rep(\cl L.,\cl H.)$ has the defined topology. That is,
$\cl L.$ is isomorphic to the algebra of continuous admissible operator
fields on $\rep(\cl L.,\cl H.).$
Using this, it is now easy to prove:
\thrm Theorem 1.5."
Let $G,\;\sigma$ be as above and let
 $\cl R.\subset\rsd$ be given. If $(G,\,\cl R.)$ has a 
$\sigma\hbox{--group}$ algebra $\cl L.,$ then up to isomorphism
it is unique."
Define a map $\wt\theta:{\cl A.(\cl L.,\cl H.)}\to
{\cl A.(\cl R.,\cl H.)}$ by $\wt\theta(T):=T\circ\theta^{-1}\,.$
That $\wt\theta$ takes admissible operator fields to 
admissible operator fields follows from Corollary~1.3.
Since $\theta$ is bijective, $\wt\theta$ is a
*-isomorphism of C*-algebras. Since $\theta$ is a
homeomorphism, it maps the strong operator continuous fields 
to the strong operator continuous fields on $\cl R.,$ i.e.
$$\leqalignno{
\wt\theta(\wt{\cl L.})=\set{T\in\cl A.(\cl R.,\cl H.)},
{T:\cl R.\to\cl B.(\cl H.)\quad\hbox{is strong--operator continuous}}..\cr}$$
But now since we have defined $\cl L.\cong\wt\theta(\wt{\cl L.})$
intrinsically on $G,$ i.e. involving only $G$ and $\cl R.,$
it follows that all $\sigma\hbox{--group}$algebras $\cl L.$
are isomorphic.

\itemitem{\bf Remark.}(1) This uniqueness theorem for group algebras
generalises previous uniqueness theorems for
(twisted) group algebras of locally compact
groups e.g. the one by Packer and Raeburn~[PR].
\itemitem{(2)} Note that the proof above provides a method for constructing 
a group algebra, i.e. {\it if we know} that $(G,\,\cl R.)$ has a
$\sigma\hbox{--group}$ algebra, then we can
construct it as the set of (strong operator) continuous admissible
operator fields on $\cl R..$ Then we obtain a $\sigma\hbox{--embedding}$ 
of $G$ in the multiplier algebra of this algebra through pointwise 
multiplication of the operator fields $T_g(\pi):=\pi(g),$ $\pi\in\cl R.,$
$g\in G\,.$
\itemitem{(3)} One can conjecture an existence theorem; e.g.
if $\cl R.\subset\rep_\sigma(G_d,\cl H.)$ is closed under direct sums,
subrepresentations,
and the  equivalence in Cor.~1.3(3), and if $\cl R.$
has ``enough'' irreducible representations (e.g. each $\pi\in\cl R.$
can be written as a direct integral of irreducible representations in $\cl R.)$
then the  set of (strong operator) continuous admissible
operator fields on $\cl R.$ is a $\sigma\hbox{--group}$ algebra
$\cl L.$ for $(G,\,\cl R.).$
An encouraging fact for the proof of this, 
is that the embedding $\varphi:G\to M(\cl L.)$ is easy to obtain
because the operator fields $T_g$ preserve the continuous admissible
fields under multiplication. 

\noindent
The difficult part of the problem is of course the existence question for 
group algebras, in particular to characterize those topological groups
for which $\cl R.=\rsg$ has a $\sigma\hbox{--group}$algebra.
It is also of interest to find subsets $\cl R.\subset\rsd$ for which
$(G,\,\cl R.)$ has a $\sigma\hbox{--group}$algebra.
Apart from the known case ${(G,\,\rsg)}$ with $G$ locally compact, 
here is a class of easy examples
of such pairs.

\item{\bf Exmp.}
Let $G$ be a nonabelian topological group, and let $\pi$ be an
irreducible representation of $G$ on a Hilbert space $\cl H._\pi$
of dimension higher than one. So
$\pi:G\to UM(\cl K.(\cl H._\pi))$ since $\cl B.(\cl H._\pi)=
M(\cl K.(\cl H._\pi))$ where $\cl K.(\cl H._\pi)$ denotes the compact 
operators on $\cl H._\pi\,.$
Recall that the strict topology of $M(\cl K.(\cl H._\pi))$
coincides with the strong operator topology. Thus the *--algebra
$\sp\pi(G)$ is strictly dense (by irreducibility) in
$M(\cl K.(\cl H._\pi)).$ Let $\cl R._\pi:=\theta(\rep\cl K.(\cl H._\pi))
\subset\rep G,$ then it is obvious that $\cl K.(\cl H._\pi)$
is a group algebra for $(G,\,\cl R._\pi),$ and that
$\cl R._\pi$ is isomorphic to the set of normal representations of 
$\cl B.(\cl H._\pi)\,.$
In the case that $\pi$ is (strong operator) continuous, all the elements of
$\cl R._\pi$ will also be continuous because they are restrictions of strictly
continuous representations of $M(\cl K.(\cl H._\pi)),$ and $\pi=\varphi.$

\noindent
Whilst group algebras may not exist for a given pair $(G,\,\cl R.),$
unitary embeddings into multiplier algebras $\varphi:G\to UM(\cl L.)$
are not hard to find, as the previous example and remark~(3) demonstrate.
Given such an embedding, here is a construction by which one
can obtain related group algebras.
Let $\cl N.$ be the strict closure of $\sp\varphi(G)$ in $M(\cl L.).$
Define $\cl S.$ to be all $N\in\cl N.$ such that
${\big\|(B_\lambda-B)N\big\|}+{\big\|N(B_\lambda-B)\big\|}
\to 0$ for all strictly convergent nets
$\{B_\lambda\}\subset\sp\varphi(G),$ $B_\lambda\to B\in\cl N..$
Clearly $\cl L.\cap\cl N.\subset\cl S..$
Then $\cl S.$ is a C*-algebra, and since products of strictly convergent nets with
fixed elements of $M(\cl L.)$ are strictly convergent, by the definition we have that
$$\varphi(G)\cl S.\subset\cl S.\supset\cl S.\varphi(G)\;.$$
If $S\in\cl S.$ then by definition there is a net $\{B_\lambda\}\subset
\sp\varphi(G)$ such that
${\big\|(B_\lambda-S)N\big\|}+{\big\|N(B_\lambda-S)\big\|}
\to 0$ for all $N\in\cl S.,$ hence there is a homomorphism 
$\psi:\sp\varphi(G)\to M(\cl S.)$ such that
${\psi\big(\sp\varphi(G)\big)}$ is strictly dense in $M(\cl S.).$
Let $\wt\theta$ denote the extension map of $\rep\cl S.$
from $\cl S.$ to $\varphi(G)\,,$
and $\wt\theta(\rep\cl S.)=:\wt{\cl R.}\subset\rsd\;.$
Then $\cl S.$ is a $\sigma\hbox{--group}$algebra
for the pair $(G,\,\wt{\cl R.}).$

Group algebras do not behave naturally w.r.t. containment,
i.e. if $\cl L._i$ is a $\sigma\hbox{--group}$algebra for
$(G,\,{\cl R._i}),$ $i=1,\,2$ where $\cl R._1\subset\cl R._2,$
then it does not always follow that $\cl L._1\subset\cl L._2$
with $\varphi_2(G)\rest\cl L._1=\varphi_1(G).$
This is because:
\thrm Proposition 1.6." Let $\cl L._i$ be a \sga
 for $(G,\,{\cl R._i}),$ $i=1,\,2$ such that
$\cl L._1\subset\cl L._2,$ and such that
$\varphi_1(g)A=\varphi_2(g)A$ for all $g\in G,\ab\;A\in\cl L._1.$
Then $\cl L._1$ is a closed two-sided ideal
of $\cl L._2,$ and hence  $\rep\cl L._2=\rep\cl L._1\oplus\rep\left(\cl L._2
\big/\cl L._1\right)$ where $\rep\cl L._1$ is identified in $\rep\cl L._2$
by unique extensions, and $\rep\left(\cl L._2
\big/\cl L._1\right)$ corresponds to those representations which vanish
on $\cl L._1.$"
Recall that $M(\cl L._1)\subset\cl L._1''\subset\cl L._2''\supset M(\cl L._2)\,.$
Recall that $\sp\varphi_i(G)$ is $\cl L._i\hbox{--strictly}$dense in $M(\cl L._i).$
Since the actions of both $\varphi_i(G),\; i=1,\,2$ coincide on
$\cl L._1,$ it follows that $\cl A.:=\sp\varphi_2(G)$ is 
$\cl L._i\hbox{--strictly}$dense in $M(\cl L._i),\; i=1,\,2$
by Prop.~1.1. Since $\cl L._1\subset\cl L._2$ it now follows from the
definition of strict topologies that the $\cl L._1\hbox{--strict}$closure
of $\cl A.$ contains the $\cl L._2\hbox{--strict}$closure
of $\cl A..$ Thus $M(\cl L._1)\supseteq M(\cl L._2)\supset\cl L._2,$
and hence $\cl L._1$ is an ideal of $\cl L._2.$
The direct sum decomposition of $\rep\cl L._2$ follows from
the ideal property, cf.~[Di].

Thus we can have natural containment of group algebras only
for direct summands.

\beginsection 2. A basic function space.

Our aim in the rest of this paper is to develop a ``universal
convolution algebra'' in which we are guaranteed to find a \sga
if it exists,
for $(G,\,\cl R.),\ab\;\cl R.\subset\rsd,$ with $G$ non-locally compact.
It then makes sense to study suitable subalgebras of it
to analyze the existence question.
In Sect.~4 we will consider in the context of these convolution
algebras the existence of $\sigma\hbox{--group}$algebras for
$(G,\,\rsg),$ i.e. the classical
Gelfand--Raikov question for continuous group representation theory.


Let $\csd$ denote the $\sigma\hbox{--twisted}$ discrete group algebra, i.e. the
C*--algebra generated by unitaries $\set\delta_x,x\in G.$ such that 
$\delta_x\cdot\delta_y=\sigma(x,\, y)\,\delta_{xy}.$
There is a bijection between the nondegenerate representations of
$\csd$ and the unitary $\sigma\hbox{--representations}$ of $G,$ and it
is given by  $\wt\pi(x):=\pi(\delta_x),$ $x\in G$ for $\pi\in\rep\csd.$
 Let $R_0\subset\rep\csd$ denote the subset in bijection with 
 $\rsg.$ 


When the group $G$ is locally compact and nondiscrete, 
we have the Haar measure $\mu,$ a notion of how a
nonzero continuous
function goes to zero at infinity, hence the function space $C_0(G),$
and the Riesz--Markov theorem which identifies the dual space
$C_0(G)^*$ with regular Borel measures on $G.$
There is then a decomposition of spaces (cf. [HR1])
$$C_0(G)^*=M_c(G)\oplus M_s(G)\oplus M_d(G)$$
where $M_d(G)$ (resp. $M_c(G),\;M_s(G))$ denotes the space of 
discrete measures (resp. continuous measures absolutely continuous w.r.t. $\mu,$
continuous measures singular w.r.t. $\mu).$
Then $C_0(G)^*$ is endowed with $\sigma\hbox{--convolution}$ and involution, w.r.t.
which $M_c(G)\cong L^1(G)$ is a *--Banach subalgebra, and its C*--envelope is
the usual group algebra $\csg,$ which is nonunital.
Then $\csg$ contains $\csd(=$C*--envelope of $M_d(G))$ in its multiplier algebra
i.e. $\csd
\subset M\big(\csg\big),$ and via the unique extension of a representation 
from $\csg$ to $M\big(\csg\big)$ we obtain a bijection between 
${\rm Rep}\,\csg$ and $R_0,$ hence $\rsg.$

Another important algebra which
a locally compact $G$ has associated to it, is its Fourier--Stieltjes algebra $B(G)$
which is a complete invariant for $G,$
cf.~[Wa].    Specifically, $B(G)$ is the space of finite spans of the continuous
positive definite functions on $G,$ (or equivalently the set of coefficient functions
for all the continuous unitary representations of $G).$ 
Now $B(G)$ is a commutative algebra w.r.t. pointwise multiplication,
and as $B(G)$ has a canonical identification with the dual space
$C^*(G)^*,$ it is a Banach space and in fact a Banach algebra.
Indeed, by this canonical identification of $B(G)$ we know
that we can identify the group algebra (as a space) with a subspace of 
$B(G)^*,$ and hence this will be a good place to look for
our generalised group algebra when $G$ is not locally compact.
We will below endow $B(G)^*$ with a convolution product
and see that the inclusion of the group algebra in $B(G)^*$
is also an inclusion of  algebras.
More concretely, note that $B(G)\subset L^\infty(G)$ is in general not
complete w.r.t. the supremum norm over $G.$ In fact, by the proof of
Corollary~13.6.5~[Di], the uniform closure of the set of functions
of compact support in $B(G)$
contains all of $C_c(G),$ and hence the uniform closure of $B(G)$
(denoted by $K(G)$ henceforth)
contains $C_0(G).$ Realize $L^1(G)\subset L^\infty(G)^*$ by
$\omega_f(h):=\int f(x)h(x)\,d\mu(x),$ $h\in L^\infty(G),$ $f\in
L^1(G),$ then  $\omega_f$ is uniquely determined by its restriction
to $C_0(G)\subset K(G).$ Thus we can identify $L^1(G)$ with a subspace 
of $K(G)^*,$ a fact which we will exploit below for more general $G.$

In the case when $G$ is not locally compact, we lose the Haar measure
$\mu,$ the space $C_0(G)=\{0\},$ and the Riesz--Markov theorem does not apply.
If $G$ is totally regular (hence $C_b(G)\cong C(\beta G)$ with $\beta G$ the
Stone--\v{C}ech compactification of $G),$ it is possible to define 
convolution for the functionals $C_b(G)^*=C(\beta G)^*$ (hence a ``generalised'' group algebra,
cf. [Do]).
However on $G$ these functionals are only finitely additive, hence
correspond to charges, not measures, cf. Alexandroff theorem
[Wh]. Thus there is no Fubini theorem for these charges, and so
 there are two inequivalent
convolutions, cf.~[Py], which do not intertwine correctly 
with the natural
involution to produce a *--algebra structure. 
Furthermore, due to the absence of the Haar measure $\mu,$
there is no way to select an analogue of $M_c(G)$ amongst these functionals.
(There are alternative characterizations of $L^1(G)=M_c(G)$
which do not use the Haar measure, cf.~[DvR, Dz, Gre, BB],
but all need some condition on measures involving
compact sets, so are not particularly
useful for us).
However $B(G)$ and $K(G)$ still makes sense even when $G$ is not locally
compact, and so, given the discussion in the previous paragraph,
we will below consider $K(G)^*$ as the appropriate universe in which
to locate generalised group algebras.

For the rest of this section, $G$ need
not be locally compact. 
Define first left and right $\sigma\hbox{--translations:}$
$$(\lambda_xf)(y):=\sigma(x,\, y)\, f(xy)=:(\rho_yf)(x)$$
for a function $f:G\to\C$ and
the involution 
$f^*(x):=\overline{f(x^{-1})}.$ For a 
 given set $\cl R.\subset\rsd,$ which is closed w.r.t. direct sums, define
the set of coefficient functions:
$$B_\sigma(\cl R.):=\set f\in C_b(G_d), {f(x)=(\psi,\,\pi(x)\varphi),\;
\pi\in\cl R.;\;\psi,\,\varphi\in{\cl H._\pi}}.$$
which is clearly a linear space, and if $\cl R.\subseteq\rsg$
these functions are also continuous.
There are two natural norms on $B_\sigma(\cl R.),$ the uniform norm
$\|f\|_\infty=\sup\limits_{x\in G}|f(x)|$ and the norm on the dual space
of $\csd,$ where we identify $B_\sigma(\cl R.)$ with a subspace
of $\csd^*$ by the bijection between $\rep\csd$ and $\rsd,$ 
$\wt\pi(x):=\pi(\delta_x)$ mentioned above, i.e.
$f(x)=(\psi,\,\wt\pi(x)\varphi)$
is identified with $\check{f}(A)=(\psi,\,\pi(A)\varphi).$
Denote the latter norm by $\|\cdot\|_*$ then it is obvious that
$\|f\|_\infty\leq\|f\|_*=\sup\set|\check{f}(A)|,{A\in\csd,\;\|A\|\leq 1}.$
because all $\delta_x$ are in the unit ball of $\csd.$
Let
 $K_\sigma(\cl R.)$ (resp. $J_\sigma(\cl R.)$) denote the completion of $B_\sigma(\cl R.)$
in the norm $\|\cdot\|_\infty$ (resp. $\|\cdot\|_*),$
 then clearly $J_\sigma(\cl R.)\subseteq K_\sigma(\cl R.).$
When the context makes clear the set $\cl R.$ under investigation, we
will simplify the notation to $B_\sigma,\;K_\sigma,$ and $J_\sigma$ and when
$\sigma=1$ we will omit the subscript. We will also not distinguish
between $\pi$ and $\wt\pi.$
We collect a few easy facts.
\thrm Theorem 2.1." Let $\cl R.\subset\rsd$ be closed w.r.t. direct sums, then\chop
(\rn1) For a fixed $x,$ the maps
$\lambda_x,\;\rho_x$ and the involution preserve $B_\sigma,\ab\;J_\sigma$ and
$K_\sigma,$ and are isometries w.r.t. both the norms $\|\cdot\|_\infty$ and
$\|\cdot\|_*$.\chop
(\rn2) $\lambda_x\lambda_y=\sigma(y,\, x)\,\lambda_{yx},$ 
$\rho_x\rho_y=\sigma(x,\,y)\,\rho_{xy},$ ${}^*$ is indeed an involution
on $C_b(G_d)$
and $(\lambda_xf)^*=\rho\s{x^{-1}}.f^*.$\chop
(\rn3) With respect to pointwise multiplication we have\chop
$X_{\sigma_1}(\cl R._1)\cdot X_{\sigma_2}(\cl R._2)\subseteq X\s\sigma_1\sigma_2.
(\cl R._1\otimes\cl R._2),$
where $X$ can be either of $B_\sigma$ or $K_\sigma,$
and if $\cl R.\subset\rep G_d$ is closed w.r.t tensor products, then
$K(\cl R.)$ is a Banach *-algebra.\chop
(\rn4) If $\cl R.$ has a faithful representation, then $B_\sigma$ separates $G,$\chop
(\rn5)  every $f\in B_\sigma(\rsg)$ is
left and right $\|\cdot\|_*\hbox{-continuous,}$\ab{(hence} uniformly continuous)
i.e. $\lim\limits_{x\to e}{\|\lambda_xf-f\|_*}=0
=\lim\limits_{x\to e}{\|\rho_xf-f\|_*}$ so we obtain the appropriate norm continuity
for functions in $J_\sigma$ and $K_\sigma\,.$ \chop
(\rn6) If $G$ is locally compact, then
$C_0(G)\subset K_\sigma(\rsg).$\chop
(\rn7) If $(G,\,\cl R.)$ has a \sga  then $B_\sigma=J_\sigma\,.$"
(\rn1) 
It is obvious that
$\|\lambda_xf\|_\infty=\|f\|_\infty=\|\rho_xf\|_\infty=\|f^*\|_\infty$ for bounded $f.$
 For a coefficient function
$f(y)={(\psi,\pi(y)\varphi)},$ $\pi\in\cl R.$ we have
$$\eqalignno{(\lambda_xf)(y)&= \sigma(x,\, y)\, f(xy)=\sigma(x,\, y)(\psi,\pi(xy)\varphi)
=(\psi,\,\pi(x)\pi(y)\varphi)\cr
&=(\pi(x)^*\psi,\,\pi(y)\varphi),&(*)\cr}$$
hence $\lambda_xf\in B_\sigma.$ 
Likewise $(\rho_xf)(y)=(\psi,\,\pi(y)\pi(x)\varphi),$ hence $\rho_xf\in B_\sigma.$
Furthermore
$$f^*(x)=\overline{f(x^{-1})}=\overline{(\psi,\,\pi(x^{-1})\varphi)}
=(\varphi,\,\pi(x^{-1})^*\psi)=(\varphi,\,\pi(x)\psi),$$
hence $f^*\in B_\sigma$ and hence ${}^*$ preserves $B_\sigma$ as well.
Now from equation~(*) above, $\lambda_xf$ 
corresponds to the functional $h(A)={(\psi,\,\pi(x)\pi(A)\varphi)}=
{(\psi,\,\pi(\delta_xA)\varphi)}$ on $\csd,$ hence
since multiplication by $\delta_x$ maps the unit ball onto the unit ball,
it follows that $\|\lambda_xf\|_*=\sup\set{|h(A)|},{A\in\csd,\,\|A\|\leq 1}.
=\|f\|_*.$ Likewise we get that $\|\rho_xf\|_*=\|f\|_*=\|f^*\|_*$.
It now follows from continuity of these maps w.r.t. the two norms
that they also preserve the two closures $J_\sigma$ and $K_\sigma.$
\chop
(\rn2) $(\lambda_x\lambda_yf)(z)=\sigma(x,\,z)\,(\lambda_yf)(xz)
=\sigma(x,\,z)\,\sigma(y,\, xz)\, f(yxz)=\sigma(y,\,x)\,\sigma(yx,\, z)\,
f(yxz)=\sigma(y,\,x)(\lambda_{yx}f)(z)$ and similarly for $\rho.$
That $f\to f^*$ is antilinear and satisfies $f^{**}=f,$ is obvious.
Moreover
$$\leqalignno{(\lambda_xf)^*(y)&=\ol(\lambda_x f)(y^{-1}).
=\ol\sigma(x,y^{-1})\,f(xy^{-1}).\cr
&=\sigma(y,x^{-1})\,\ol f((yx^{-1})^{-1}).
=\sigma(y,x^{-1})\, f^*(yx^{-1})\cr
&=\rho\s x^{-1}.f^*(y)\;.}$$
(\rn3)  By continuity of the pointwise product w.r.t. the uniform norm,
it suffices to prove that:
$B_{\sigma_1}(\cl R._1)\cdot B_{\sigma_2}(\cl R._2)\subseteq B\s\sigma_1\sigma_2.(\cl R._1\otimes\cl R._2).$
Let $f_i\in B_{\sigma_i}(\cl R._i),\;\ab i=1,\,2,$ i.e.
$f_i(x)={\big(\varphi_i,\,\pi_i(x)\psi_i\big)}$ where $\pi_i\in\cl R._i,$   
$\varphi_i,\;\ab\psi_i\in\cl H.\s\pi_i..$ Then
$$f_1(x)f_2(x)=\big(\varphi_1\otimes\varphi_2,\,(\pi_1\otimes\pi_2)(x)\,\psi_1\otimes\psi_2
\big)$$
and since $\pi_1\otimes\pi_2\in \cl R._1\otimes\cl R._2,$  
it follows that 
$f_1\cdot f_2\in B\s\sigma_1\sigma_2.(\cl R._1\otimes\cl R._2).$
If $\cl R.$ is closed under tensor products (as well as direct sums), then
$K_1(\cl R.)$ is closed under pointwise multiplication, the involution, and is uniformly 
closed, so it follows that it is a Banach *-algebra.\chop
(\rn4) If there is a faithful representation
$\pi\in\cl R.,$ its coefficient functions $f(x)={(\psi,\,\pi(x)\varphi)}$ must
separate $G$ because they determine $\pi(x).$\chop
(\rn5) 
It is only necessary to establish continuity w.r.t. $\|\cdot\|_*$
for the elements of $B_\sigma(\rsg).$ For left continuity
of $f(x)={(\psi,\,\pi(x)\varphi)}:$
$$\leqalignno{\left\|\lambda_xf-f\right\|_*&
=\sup\set{|\big(\psi,\big(\pi(x)-\un\big)\pi(A)\varphi\big)|},{A\in\csd,\;\|A\|\leq 1}.\cr
&\leq\set{\big\|(\pi(x^{-1})-\un)\psi\big\|\cdot\|\pi(A)\varphi\|},{A\in\csd,\;\|A\|\leq 1}.\cr
&\leq\left\|\big(\pi(x^{-1})-\un\big)\psi\right\|
\cdot\|\varphi\|\maprightt x,e.0\cr}$$
from which it is clear, and similarly we get right continuity.\chop
(\rn6) Let $G$ be locally compact, then
construct the usual group extension $G_\sigma$ and recall the 
bijection between $\rsg$ and ${\rm Rep}\,G_\sigma$ (cf.~[Ma]),
given by $\pi_\sigma(x,\,t):=\pi(x)t,$ $\pi\in\rsg,$ $x\in G,$ $t\in{\Bbb T}.$
From the argument at the start of this section
(using the proof of Corollary~13.6.5 in ~[Di]),
 we have that $C_0(G_\sigma)\subset K(G_\sigma).$
Now since for a coefficient function on $G_\sigma$ we have
${(\psi,\,\pi_\sigma(x,\,t)\varphi)}={t(\psi,\,\pi(x)\varphi)},$
it is clear that the restriction of $K(\rep G_\sigma)$ to $G\subset G_\sigma$ is
just $K_\sigma(\rsg).$ Thus restriction to $G$ produces
$C_0(G)\subset K_\sigma(\rsg),$ as required.\chop
(\rn7) Since 
 $(G,\,\cl R.)$ has a \sga $\cl L.,$ then  by Proposition~1.1
each $f\in B_\sigma$ has a unique strictly continuous extension 
$\wh{f}$ from $\sp\varphi(G)$ to $M(\cl L.).$ Now
$$\|f\|_*=\left\|\wh{f}\restriction\varphi(\csd)\right\|=\|\wh{f}\|=
\|\wh{f}\restriction\cl L.\|$$
because $\wh{f}$ is strictly continuous, both $\varphi(\csd)$ and $\cl L.$
are strictly dense in $M(\cl L.)$ and the unit ball of any strictly dense C*-algebra
in $M(\cl L.)$ is strictly dense in the unit ball of $M(\cl L.)$
(the last fact is Exercise~2.N in [WO]).
But $\cl L.$ is a group algebra for $\cl R.,$ hence
$$\set{\wh{f}\restriction\cl L.},f\in B_\sigma.=\cl L.^*$$
and this is complete in norm. Thus $B_\sigma$ is complete in the
${\|\cdot\|_*}\hbox{--norm}$ and hence $B_\sigma=J_\sigma\,.$

\noindent In the rest of this paper we will always assume that the sets
$\cl R.\subseteq\rsd$ under consideration (for construction of group algebras)
are closed with respect to direct sums.

\beginsection 3. Convolution Algebras.

Fix a set $\cl R.\subseteq\rsd$ (closed with respect to direct sums).
Here we want to make the dual spaces $J_\sigma(\cl R.)^*$ and $K_\sigma(\cl R.)^*$ into
convolution algebras, following the method of Def.~19.1~[HR1].
Since we can identify $J_\sigma(\cl R.)^*$ (resp. $K_\sigma(\cl R.)^*)$
with the functionals on $B_\sigma(\cl R.)$ which are
$\|\cdot\|_*\hbox{--continuous}$\ab(resp. $\|\cdot\|_\infty\hbox{--continuous}),$
it follows from $\|\cdot\|_\infty\leq\|\cdot\|_*$ that
$K_\sigma^*\subseteq J_\sigma^*\,.$
Before considering the convolutions, we need notation and a
preparatory lemma.
\chop {\bf Notation:} Let $f\in C_b(G\times\cdots\times G)$ (n factors)
and let $L\subseteq C_b(G)$ be a closed subspace. Assume that the
function $x_1\to f(x_1,\ldots,\, x_n)$ for fixed $x_j,\;j\not=1$
is in $L.$ Now, given functionals $\omega^{(1)},\ldots,\,\omega^{(n)}\in L^*,$
the expression:
$$\omega^{(n)}\s x_n.\left(\omega^{(n-1)}\s x_{n-1}.\left(
\cdots\left(\omega^{(1)}\s x_1.\big(f(x_1,\ldots,\,x_n)\big)\right)\cdots\right)
\right)$$
means the following. Starting with the innermost functional
$\omega^{(1)},$ first evaluate $\omega^{(1)}$ of the function
$x_1\to f(x_1,\,x_2,\ldots,\,x_n)$ obtaining a function 
$F(x_2,\ldots,\,x_n).$ Next evaluate $\omega^{(2)}$ of the function
$x_2\to F(x_2,\ldots,\,x_n)$ to get $G(x_3,\ldots,\,x_n).$
Continue until all functionals are evaluated. 
For this to make sense, we need to have that all subsequent 
 functions $x_2\to F(x_2,\ldots,x_n)$ etc., to be evaluated are in $L.$
This brings the functional notation closer to integral notation,
e.g. $$\omega_x(f(x))=\omega(f)=\int f(x)\,d\mu(x)$$
for a functional $\omega$ given by a measure.
\thrm Lemma 3.1." For each $\omega\in J_\sigma^*$ and $\pi\in\cl R.,$
there is a unique operator $\pi(\omega)\in\cl B.(\cl H._\pi)$
such that $\|\pi(\omega)\|\leq\|\omega\|$ and
${\omega_x\big((\psi,\,\pi(x)\,\varphi)\big)}={(\psi,\,\pi(\omega)\,\varphi)}$
for all $\psi,\;\varphi\in\cl H._\pi.$
Moreover $\pi(\omega)$ preserves each cyclic component
of $G$ in $\cl H._\pi,$ i.e. $\pi(\omega)\cl H._\psi
\subseteq\cl H._\psi$ for all $\psi\in\cl H._\pi,$
where $\cl H._\psi:=[\pi(G)\psi].$ "
Let $\omega\in J_\sigma^*,$ then 
${\omega_x\big((\psi,\,\pi(x)\,\varphi)\big)}$ exists since
$x\to{(\psi,\,\pi(x)\,\varphi)}$ is in $B_\sigma.$
Now the map $\psi\to{\omega_x\big((\psi,\,\pi(x)\,\varphi)\big)}$
is conjugate linear, and bounded as
$$\eqalignno{\left|{\omega_x\big((\psi,\,\pi(x)\,\varphi)\big)}\right|&\leq
\|\omega\|\cdot\sup\set{\left|(\psi,\,\pi(A)\,\varphi)\right|},{A\in\csd,\,
\|A\|\leq 1}.\cr
&\leq
\|\omega\|\cdot\|\psi\|\cdot\|\varphi\|&{(*)}\cr}$$
hence it is a conjugate linear functional on $\cl H._\pi.$
Thus by the Riesz representation theorem, there is a vector
$\varphi_\omega\in\cl H._\pi$ such that
$${\omega_x\big((\psi,\,\pi(x)\,\varphi)\big)}=(\psi,\,\varphi_\omega)\quad
\forall\;\psi\in\cl H._\pi\eqno{(+)}$$
Denote $\varphi_\omega$ by $\pi(\omega)\varphi,$ then by $(*)$ we see
$\|\pi(\omega)\varphi\|\leq\|\omega\|\cdot\|\varphi\|,$
hence by linearity of $\varphi\to\pi(\omega)\varphi$ (clear from $(+)),$
we have defined a bounded operator $\pi(\omega):\cl H._\pi\to\cl H._\pi.$
Uniqueness comes from the fact that $\pi(\omega)$ is fully determined by
the coefficients $(\psi,\,\pi(\omega)\varphi)$ as $\psi$ and
$\varphi$ ranges over $\cl H._\pi.$\chop
Finally, fix $\psi\in\cl H._\pi\backslash 0,$ then the functional
$\xi\to{\omega_x\big((\xi,\,\pi(x)\,\varphi)\big)}$
for $\varphi\in\cl H._\psi$ is zero on the orthogonal 
complement of $\cl H._\psi,$ hence the vector $\varphi_\omega
=\pi(\omega)\varphi$ from the Riesz theorem must be 
in $\cl H._\psi.$

Note that when $\omega$ is associated to a measure: $\omega(f)=
\int f\,d\mu,$ then
$$\pi(\omega)=\int\pi(x)\,d\mu(x)\;.$$

In the case that $G$ is locally compact, the convolution of two measures 
$\mu,\;\nu\in M(G)=C_0(G)^*$ is given by:
$$\eqalignno{\int f(x)\, d(\mu*\nu)(x)&=
\int\int\sigma(x,\, y)\,f(xy)\,d\mu(x)\,d\nu(y) &(1)\cr
&=\int\int\sigma(x,\, y)\,f(xy)\,d\nu(y)\,d\mu(x)\quad\forall\; f\in
C_0(G), &(2)\cr}$$
by Fubini. (When we generalise to non-locally compact groups 
or other types of functionals, 
these two formulii  will give different convolutions~[Py].)
Let us first rewrite these in terms of the associated functionals
$\omega(f):=\int f d\mu$ and $\xi(f):=\int f\,d\nu.$
For $f\in C_0(G),$ $\psi\in C_0(G)^*,$ define
$$\eqalignno{f^\psi(x):=\psi(\lambda_xf)=\psi_y(\sigma(x,y)f(xy))\qquad\hbox{and}
\qquad f_\psi(x):=\psi(\rho_xf)=\psi_y(\sigma(y,x)f(yx))}$$
Then by Lemma~19.5~[HR1], $f^\psi\in C_0(G)\ni f_\psi,$
hence we can write the convolution formulii as
$$\eqalignno{(\omega*\xi)(f)&:= \int f(x)\,d(\mu*\nu)(x)
=\int\int(\rho_yf)(x)\,d\mu(x)\,d\nu(y)\cr
&=\xi(f_\omega) = \xi_y(\omega_x(\sigma(x,y)\,f(xy)))&(3)\cr
\hbox{and}\qquad\qquad\quad\cr
(\omega*\xi)(f)&= \int\int(\lambda_xf)(y)\,d\nu(y)\,d\mu(x)\cr
&=\omega(f^\xi)=\omega_x(\xi_y(\sigma(x,y)\,f(xy)))&(4)\cr}$$
Henceforth, let $G$ be non--locally compact. In order to generalise
(3) and (4) to $G,$ we first need the lemma:
\thrm Lemma 3.2." Let $X$ denote either $J$ or $K.$ 
Let $f\in X_\sigma(\cl R.),$ $\omega\in X_\sigma^*$ and define
as above $f^\omega(x):=\omega(\lambda_xf)$ and $f_\omega(x):=\omega(\rho_xf).$
Then $f^\omega\in X_\sigma(\cl R.)\ni f_\omega.$"
Consider $f\in B_\sigma$ of the form $f(x)=(\psi,\,\pi(x)\varphi),$ $\pi\in\cl R..$
Let $\omega\in X_\sigma^*$ then
$$\leqalignno{\qquad
f^\omega(x)&=\omega(\lambda_xf)=\omega_y\left(\sigma(x,\, y)\,
\big(\psi,\,\pi(xy)\varphi\big)\right)\cr
&=\omega_y\big((\psi,\,\pi(x)\pi(y)\varphi)\big)=\omega_y\big((\pi(x)^*\psi,\,
\pi(y)\varphi)\big)\cr
&=\big(\pi(x)^*\psi,\,\pi(\omega)\varphi\big)=\big(\psi,\,\pi(x)\pi(\omega)\varphi\big)
\cr}$$
making use of Lemma~3.1. Thus $f^\omega\in B_\sigma.$ For fixed $\omega$
the map $f\to f^\omega$ is linear and norm--continuous by
$$\|f^\omega\|_\infty=\sup_x\big|\omega(\lambda_xf)\big|\leq
\|\omega\|\cdot\sup_x\|\lambda_xf\|_\infty=\|\omega\|\cdot\|f\|_\infty$$
in the case $X=K.$ For the case  $X=J,$ consider an $f\in B_\sigma$ as above,
then 
$$\eqalignno{\|f^\omega\|_*&=\sup\set{|\omega_y\big((\psi,\,\pi(A)\pi(y)\varphi)\big)|},
{A\in\csd,\,\|A\|\leq 1}.\cr
&\leq\sup\set{\|\omega\|\cdot\|h_A\|_*},{A\in\csd,\,\|A\|\leq 1}.\cr}$$
where $h_A(y):=(\psi,\,\pi(A)\pi(y)\varphi).$
Now 
$$\|h_A\|_*=\sup\set{|(\psi,\,\pi(A)\pi(B)\varphi)|},
{B\in\csd,\,\|B\|\leq 1}.\leq\|f\|_*$$
because $\|AB\|\leq\|A\|\|B\|\leq 1.$ 
Thus $\|f^\omega\|_*\leq\|\omega\|\|f\|_*,$ and hence we get that
$f^\omega\in X_\sigma$ for all $f\in X_\sigma\,.$
Likewise we have:
$$\eqalignno{f_\omega(x)&=\omega_y((\rho_xf)(y))=
\omega_y\big((\psi,\,\pi(y)\pi(x)\varphi)\big)\cr
&=(\psi,\,\pi(\omega)\pi(x)\varphi)=(\pi(\omega)^*\psi,\,\pi(x)\varphi)\cr}$$
from which it is clear that $f_\omega\in B_\sigma,$ and hence
by a similar argument as above, $f\to f_\omega$ is norm continuous, so
$f_\omega\in X_\sigma$ for all $f\in X_\sigma.$

Now as remarked before,
there are two possible convolutions one can define on $J_\sigma^*$
given by (3) and (4),
but surprisingly they are the same:

\thrm Theorem 3.3."(\rn1) Given $\omega,\;\xi\in J_\sigma^*,\; f\in J_\sigma,$
we have $\xi(f^\omega)= \omega(f_\xi),$ i.e.
$$\xi_x\big(\omega_y\big(\sigma(x,\, y)\, f(xy)\big)\big)
=\omega_y\big(\xi_x\big(\sigma(x,\, y)\, f(xy)\big)\big),$$
and thus we can define convolution in $J_\sigma^*$ by
$$(\xi*\omega)(f):=\xi(f^\omega)= \omega(f_\xi).$$\chop
Moreover, $K_\sigma^*\subseteq J_\sigma^*$ is closed under this convolution.\chop
(\rn2) Both $J_\sigma^*$ and $K_\sigma^*$ are Banach *--algebras w.r.t. convolution,
the involution $\omega^*(f):=\overline{\omega(f^*)},$
 and their usual norms.\chop
(\rn3) Let $\pi\in\cl R.,$ and let $\omega\to\pi(\omega)$ be the map
from $J_\sigma^*$ to $\cl B.(\cl H._\pi)$ given in Lemma~3.1.
Then this map is a continuous *--homomorphism of the Banach *-algebra
$J_\sigma^*,$ i.e. a representation and hence it restricts to
a continuous representation of $K_\sigma^*.$"
(\rn1) By Lemma 3.2, both of $\xi(f^\omega)$ and $\omega(f_\xi)$
exist. Let $f\in B_\sigma$ be of the form $f(x)=(\psi,\,\pi(x)\varphi),$
$\pi\in\cl R.,$
then 
$$\eqalignno{\xi(f^\omega)&=\xi_x\big(\omega_y\big(\sigma(x,\, y)\, f(xy)\big)\big)\cr
&=\xi_x\big(\omega_y\big((\psi,\,\pi(x)\pi(y)\varphi)\big)\big)
=\xi_x\big((\pi(x)^*\psi,\,\pi(\omega)\varphi)\big)\cr
&=(\psi,\,\pi(\xi)\pi(\omega)\varphi)=(\pi(\xi)^*\psi,\,\pi(\omega)\varphi)\cr
&=\omega_y\big((\psi,\,\pi(\xi)\pi(y)\varphi)\big)
=\omega_y\big(\xi_x\big((\psi,\,\pi(x)\pi(y)\varphi)\big)\big)\cr
&=\omega_y\big(\xi_x\big(\sigma(x,\, y)\, f(xy)\big)\big)=\omega(f_\xi).
&(5)\cr}$$   
By the norm--continuity found in the proof of Lemma~3.2:
$$\|\xi(f^\omega)\|\leq\|\xi\|\cdot\|\omega\|\cdot\|f\|_*\geq\|\omega(f_\xi)\|,$$
the equality (5) extends to all of $J_\sigma,$ thus establishing the claim.\chop
If $\omega,\,\xi\in K_\sigma^*,$ then so is $\xi*\omega,$ by the continuity
$\|\xi(f^\omega)\|\leq{\|\xi\|\cdot\|f^\omega\|_\infty}
\leq{\|\xi\|\cdot\|\omega\|\cdot\|f\|_\infty}$ encountered
above.\chop
(\rn2) 
Linearity of $\xi*\omega$ in $\xi,\,\omega\in J_\sigma^*$ is clear from the
definition. For associativity:
$$\eqalignno{\big((\xi*\beta)*\omega\big)(f)&=
(\xi*\beta)_x\big(\omega_y\big(\sigma(x,\,y)\,f(xy)\big)\big)\cr
&=\xi_z\big(\beta_v\big(\sigma(z,\,v)\omega_y\big(\sigma(zv,\,y)\, f(zvy)
\big)\big)\big)\cr
&=\xi_z\big(\beta_v\big(\omega_y\big(\sigma(v,\,y)\,\sigma(z,\,vy)\, f(zvy)
\big)\big)\big)\cr
&=\xi_z\big((\beta*\omega)_x\big(\sigma(z,\, x)\, f(zx)\big)\big)\cr
=\big(\xi*(\beta*\omega)\big)(f)\cr}$$
making use of the two--cocycle relation for $\sigma.$
Norm continuity of course follows from ${|(\xi*\omega)(f)|}\leq
\|\xi\|\|\omega\|\|f\|_*$  or
${|(\xi*\omega)(f)|}\leq
\|\xi\|\|\omega\|\|f\|_\infty$ if instead $\xi,\,\omega\in K_\sigma^*.$
Thus $J_\sigma^*$ and $K_\sigma^*$ are Banach algebras.
Concerning the involution, it is clear that $\omega\to\omega^*$ is
antilinear and that $\omega^{**}=\omega.$
Now 
$$\eqalignno{\|\omega^*\|&=\sup\set\big|\omega^*(f)\big|,
f\in K_\sigma,\;\|f\|\leq 1.\cr
&= \sup\set\big|\omega(f^*)\big|,
f\in K_\sigma,\;\|f\|\leq 1.=\|\omega\|\cr}$$
using $\|f^*\|=\|f\|$ where here $\|f\|$ denoted 
$\|f\|_\infty$ if $\omega\in K_\sigma^*,$ and
$\|f\|_*$ otherwise.
Furthermore, for $\omega,\;\xi\in J_\sigma^*,$
$$\eqalignno{(\xi*\omega)^*(f)&=\overline{(\xi*\omega)(f^*)}
=\overline{\xi_x\big(\omega_y\big(\sigma(x,\,y)(f^*)(xy)}\cr
&=\overline{\xi_x\big(\omega_y\big(\sigma(x,\,y)\overline{
f(y^{-1}x^{-1})}\big)}\cr
&=\overline{\xi_x\big(\omega_y\big(\overline{\sigma(y^{-1},x^{-1})\,
f(y^{-1}x^{-1})}\big)}\cr
&=\overline{\xi_x\big(\overline{\omega^*_y\big(\sigma(y,x^{-1})\,
f(yx^{-1})}\big)}\cr
&=\xi^*_x\big(\omega^*_y\big(\sigma(y,x)\,
f(yx)\big)=(\omega^**\xi^*)(f)\cr}$$
where we made use of $\overline{\sigma(x,\,y)}=\sigma(y^{-1},\,x^{-1})\,.$
Thus $J_\sigma^*$ and $K_\sigma^*$ are Banach *--algebras.\chop
(\rn3) That $\omega\to\pi(\omega)$ for $\pi\in\cl R.$ is linear is 
easy to see. From Lemma~3.1 we also have that $\|\pi(\omega)\|\leq\|\omega\|,$
hence the map is continuous. 
We show that $\omega\to\pi(\omega)$ is a homomorphism.
$$\eqalignno{(\psi,\,\pi(\omega*\xi)\varphi)&=
(\omega*\xi)_x\big((\psi,\,\pi(x)\varphi)\big)\cr
&=\omega_x\big(\xi_y\big(\sigma(x,\,y)\,(\psi,\,\pi(xy)\varphi)\big)\big)\cr
&=\omega_x\big(\xi_y\big((\pi(x)^*\psi,\,\pi(y)\varphi)\big)\big)\cr
&=\omega_x\big((\psi,\,\pi(x) \pi(\xi)\varphi)\big)
=(\psi,\,\pi(\omega) \pi(\xi)\varphi)\;.\cr}$$
Finally, to establish that $\pi$ is a *--homomorphism of $J_\sigma^*,$
$$\eqalignno{(\psi,\,\pi(\omega^*)\varphi)&=
(\omega^*)_x\big((\psi,\,\pi(x)\varphi)\big)\cr
&=\overline{\omega_x\big((\psi,\,\pi(x)\varphi)^*\big)}
=\overline{\omega_x\big((\pi(x^{-1})\varphi,\,\psi)\big)}\cr
&=\overline{\omega_x\big((\varphi,\,\pi(x)\psi)\big)}
=\overline{(\varphi,\,\pi(\omega)\psi)}\cr
&=(\psi,\,\pi(\omega)^*\varphi)\cr}$$
hence $\pi(\omega^*)=\pi(\omega)^*.$

There are several distinguished subalgebras of $J_\sigma^*,$ but
one in particular which we'll find useful, is the algebra of
point measures, i.e. $\delta_x\in K_\sigma^*\subseteq J_\sigma^*,\; x\in G$ defined by
$\delta_x(f):=f(x),\;\ab\;f\in K_\sigma.$

\thrm Theorem 3.4." For any $\omega\in J_\sigma^*,$ we have that
$\delta_x*\omega=\omega\circ\lambda_x$ and
${\omega*\delta_x}=\omega\circ\rho_x,$ hence $\delta_e$ is the identity
$\un$ of $J_\sigma^*.$ Moreover $\delta_x^*=\delta\s x^{-1}.=(\delta_x)^{-1},$
and $\delta_x*\delta_y=\sigma(x,\, y)\,\delta_{xy}$."
$(\delta_x*\omega)(f)=(\delta_x)_z\big(\omega_y\big(\sigma(z,\,y)\,
f(zy)\big)\big)=\omega_y\big(\sigma(x,\,y)\,f(xy)\big)
=\omega(\lambda_xf).$ Likewise $\omega*\delta_x=\omega\circ\rho_x.$\chop
$(\delta_x^*)\s y.(f(y))=\overline{(\delta_x)\s y.(\overline{f(y^{-1})})}
=f(x^{-1})=\delta\s x^{-1}.(f)$ i.e. $\delta_x^*=\delta\s x^{-1}..$
\chop Furthermore,
$$\eqalignno{(\delta_x*\delta_y)(f)&=(\delta_x)_z\big((\delta_y)_v
(\sigma(z,\, v)\,f(zv))\big)=
(\delta_x)_z\big(\sigma(z,\, y)\,f(zy)\big)\cr
&=\sigma(x,\, y)\,f(xy)=\big(\sigma(x,\, y)\,\delta\s xy.\big)(f),\cr}$$
i.e. $\delta_x*\delta_y=\sigma(x,\, y)\,\delta_{xy}$ and hence
$\delta\s x^{-1}.=(\delta_x)^{-1}=\delta_x^*.$

In the case that $K_\sigma(\cl R.)$ separates $G,$ we get that the
$\sigma\hbox{--homomorphism}$
$\delta:G\to J_\sigma^*$ is injective, and then
the C*--enveloping algebra of the *-algebra $\sp\delta\s G.$
is the (twisted) discrete
group algebra $C^*(G_d).$ 

The next theorem establishes that the convolution algebra
$J_\sigma^*$ is the appropriate setting in which to
locate group algebras, if they exist.

\thrm Theorem 3.5."
Let $\cl L.$ be a \sga for the pair $(G,\,\cl R.)$ with 
associated $\sigma\hbox{--homomorphism}$ $\varphi:G\to UM(\cl L.),$
then there is a 
*-isomorphism $\Psi:M(\cl L.)\to J_\sigma(\cl R.)^*$ (into) 
such that 
$\Psi(\varphi(x))=\delta_x$
for all $x\in G.$"
Given the group algebra $\cl L.,$ recall from Proposition~1.1 that 
each $\pi\in\cl R.$ has a unique strict--strong operator continuous
extension from $\varphi(G)$ to a representation $\wh\pi\in\rep M(\cl L.),$ 
hence each $f\in B_\sigma$ has a canonical extension to
$M(\cl L.)$ by $\wh{f}(A):={(\psi,\,\wh\pi(A)\xi)},$ $A\in\cl L.$
when $f(x)={(\psi,\,\pi(x)\xi)}.$ 
Now $\wh\pi\left(C^*(\varphi(G))\right)=\pi\left(\csd\right),$ so
 by definition $\|f\|_*={\left\|\wh{f}\rest C^*(\varphi(G))\right\|}$
so since by Prop.~1.1 $C^*(\varphi(G))$ is strictly dense in $M(\cl L.),$
it follows by the strict-strong operator continuity of $\wh\pi$ that
${\|\wh{f}\|}\leq\|h\|$ where $h$ is the functional
$h(A):={(\psi,\,(A)\xi)}$ on $\pi(G)''.$ By a Kaplansky
density argument we know that $\|h\|={\left\|h\rest
\wh\pi\big(C^*(\varphi(G))\big)\right\|}=\|f\|_*$ and thus
 ${\|\wh{f}\|}=\|f\|_*.$
Hence the map $f\to\wh{f}$ is a linear isometry,
$\Phi:J_\sigma\to M(\cl L.)^*,$ and in particular
by $M(\cl L.)\subset M(\cl L.)^{**}$ we obtain a linear map
$\Psi:M(\cl L.)\to J_\sigma^*$ by
$\Psi(A)(f):=\Phi(f)(A)$ for $f\in J_\sigma,$
$A\in M(\cl L.),$ i.e. on $B_\sigma$
$$\Psi(A)_x\big((\psi,\,\pi(x)\xi)\big)=(\psi,\,\wh\pi(A)\xi)
=\big(\psi,\,\pi(\Psi(A))\xi)$$
by Lemma 3.1, hence $\wh\pi(A)=\pi(\Psi(A))\,.$
This establishes also that $\Psi(A)$ is $\|\cdot\|_*\hbox{--continuous}$
on $B_\sigma,$ hence confirms that $\Psi(A)\in J_\sigma^*.$
In particular if $A=\varphi(x)$ and $f(x)={(\psi,\,\pi(x)\xi)}$ then
$$\Psi(A)(f)=(\psi,\,\wh\pi(\varphi(x))\xi)=(\psi,\,\pi(x)\xi)
=(\delta_x)_y\big((\psi,\,\pi(y)\xi)\big)$$
and by letting $f$ range over $B_\sigma$ we get that
$\Psi(\varphi(x))=\delta_x\,.$ 
To see that $\Psi$ is a homomorphism:
$$\eqalignno{\Psi(A)*\Psi(B)(f)&=
\left(\Psi(A)*\Psi(B)\right)_x\big((\psi,\,\pi(x)\xi)\big)
=\left(\psi,\,\pi\big(\Psi(A)*\Psi(B)\big)\xi\right)\cr
&=\left(\psi,\,\pi\big(\Psi(A)\big)\pi\big(\Psi(B)\big)\xi\right)
&\hbox{by Theorem 3.3(\rn3)}\cr
&=(\psi,\,\wh\pi(AB)\xi)=\Psi(AB)(f)\,.\cr}$$
The adjoint is preserved because
$$\eqalignno{\Psi(A)^*(f)&=\ol\Psi(A)(f^*).=
\ol\Psi(A)_x\big(\xi,\,\pi(x)\psi\big).
=\ol\big(\xi,\,\wh\pi(A)\psi\big).\cr
&=\big(\psi,\,\wh\pi(A^*)\xi\big)=\Psi(A^*)(f)\,.}$$
Thus $\Psi:M(\cl L.)\to J_\sigma^*$ is a *-homomorphism of a 
C*-algebra into a Banach *-algebra, and hence
$\|\Psi(A)\|\leq\|A\|$ for all $A\in M(\cl L.)\,.$
In fact $\Psi$ is isometric on $\cl L.$ because for all
$A\in \cl L.:$
$$\eqalignno{\big\|\Psi(A)\big\|&=
\sup\set|\Psi(A)(f)|,f\in J_\sigma,\;\|f\|_*\leq 1.\cr
&=\sup\set|\Psi(A)(f)|,f\in B_\sigma,\;\|f\|_*\leq 1.&\hbox{as $B_\sigma$ is dense in $J_\sigma.$}\cr
&=\sup\set{\left|\big(\psi,\,\wh\pi(A)\xi\big)\right|},{\pi\in\cl R.;\;\psi,\xi\in\cl H._\pi,\;
\|\psi\|\leq 1\geq\|\xi\|}.\cr
&=\sup\set{\left|\big(\psi,\,\pi(A)\xi\big)\right|},{\pi\in\rep\cl L.;\;\psi,\xi\in\cl H._\pi,\;
\|\psi\|\leq 1\geq\|\xi\|}.\cr
&=\|A\|\,\cr}$$
where in the penultimate step we used the fact that $\wh\pi\rest\cl L.=\theta^{-1}(\pi)$ and that
$\theta:\rep\cl L.\to\cl R.$ is a bijection.
Thus $\Psi$ is an isomorphism on $\cl L.,$ and as multipliers $B\in M(\cl L.)$
are uniquely determined by their action on $\cl L.,$ it follows that $\Psi$
must also be an isomorphism on $M(\cl L.).$

By this theorem, one should look for group algebras for a given pair $(G,\,\cl R.)$
in the subalgebras of $J_\sigma^*$ which are stable under multiplication
by $\delta\s G.\,.$
\item{\bf Def.} A {\it d--ideal} $\cl A.$ of $J_\sigma(\cl R.)^*$ is a
nonzero norm--closed *--subalgebra
such that $\delta_x*\cl A.\subseteq\cl A.\supseteq\cl A.*\delta_x$
for all $x\in G,$
(i.e. $\delta\s G.$ is in the relative multiplier algebra of $\cl A..)$

\noindent
From any nonzero $A\in J_\sigma^*$ we can generate a d-ideal
by just taking the closed *-algebra generated by 
the set $\delta\s G.*A.$ For any d-ideal we have the usual map
$\theta:\rep\cl A.\to \rsd$ by $\theta(\pi)(x)=\slim\limits_\alpha
\pi(\delta_x*E_\alpha)$ where $\{E_\alpha\}\subset\cl A.$ is any
approximate identity of $\cl A..$ 
 (Equivalently, $\theta(\pi)$ is uniquely determined by the 
 equation $\theta(\pi)(x)\cdot\pi(A)\psi=\pi(\delta_x*A)\psi$ 
 for all $A\in\cl A.,$ $\psi\in\cl H._\pi.)$

Note that in the proof of Theorem~3.5 we established that $\|\Psi(A)\|=\|A\|$ on $\cl L.,$
hence on the image of the group algebra in $J_\sigma^*,$ the norm is already
a C*-norm. This has a striking generalization:

\thrm Theorem 3.6."
The norm of $J_\sigma(\cl R.)^*$  satisfies
$$\|A\|=\sup\set\|\pi(A)\|,{\pi\in\cl R.}.$$ for $A\in J_\sigma(\cl R.)^*.$
 This is a C*-norm, hence $J_\sigma(\cl R.)^*$
is a C*-algebra, and so is every d-ideal in $J_\sigma(\cl R.)^*.$"
We adapt the
same calculation at the end of the proof of Theorem~3.5 to this case.
Recall first from Lemma~3.1 and Theorem~3.3(\rn3) that each ${\pi\in\cl R.}$
defines a *-representation of the Banach *-algebra
 $J_\sigma(\cl R.)^*$ and that $A(f)=
{\big(\psi,\,\pi(A)\xi\big)}$ for any coefficient function
$f(x) ={\big(\psi,\,\pi(x)\xi\big)}$ and $A\in J_\sigma^*.$
Now for all
$A\in J_\sigma^*:$
$$\eqalignno{\big\|A\big\|&=
\sup\set|A(f)|,f\in J_\sigma,\;\|f\|_*\leq 1.\cr
&=\sup\set|A(f)|,f\in B_\sigma,\;\|f\|_*\leq 1.&\hbox{as $B_\sigma$ is dense in $J_\sigma$}\cr
&=\sup\set{\left|\big(\psi,\,\pi(A)\xi\big)\right|},{\pi\in\cl R.;\;\psi,\xi\in\cl H._\pi,\;
\|\psi\|\leq 1\geq\|\xi\|}.\cr
&=\sup\set{\left\|\pi(A)\right\|},{\pi\in\cl R.}..\cr}$$
Since the operator norms $\|\pi(A)\|$ are C*-norms, it follows that $\|\cdot\|$ on 
$J_\sigma(\cl R.)^*$ is a C*-norm.

By this theorem our search for group algebras is simplified, in that the d-ideals
under consideration are semisimple (i.e. have zero radicals) and are already closed in
a C*-norm, so it is unnecessary to consider C*-envelopes.


\beginsection 4. Conditions for group algebras.

Inspired by Theorem 3.5, we now want to characterise the properties which a
{d-ideal} in $J_\sigma^*$ should satisfy in order to be a \sga
for $(G,\,\cl R.).$ We will specialise to the case $\cl R.=\rsg$
(the Gelfand--Raikov question) at the end of this section.

%
Obviously a d-ideal must satisfy
$\theta(\rep\cl A.)\subseteq\cl R.$ if it is to be a \sga for
$(G,\,\cl R.).$ So we denote
$$\cl I.(\cl R.):=\set\cl A.\subset J_\sigma^*,{\cl A.\quad\hbox{is a d-ideal and}\quad
\theta(\rep\cl A.)\subseteq\cl R.}.\,.$$
\def\pa{\pi\s{\cl A.}.}
(This set will be analyzed later for $\cl R.=\rsg).$
The natural map which we will want to be inverse to $\theta,$ is the map
${\pi\in\rsg\to\pa\in\rep\cl A.}$ defined by
$$\omega_x\big((\psi,\,\pi(x)\varphi)\big)=(\psi,\,\pa(\omega)\varphi)\qquad\quad
\forall\;\psi,\,\varphi\in\cl H._\pi,\;\omega\in\cl A.\;,$$
via Lemma 3.1 and Theorem 3.3(\rn3).
In the case that $G$ is locally compact and
$\cl R.=\rsg,$ we realize $L^1(G)$ in $J_\sigma^*$ 
as usual by integrals:
$\omega_h(f):=\int_G f(x)h(x)\,d\mu(x),\ab\;f\in B_\sigma,\ab\;h\in L^1(G)$
with $\mu$ the Haar measure
(then $\omega_h\in K_\sigma^*\subset J_\sigma^*).$ Then the closure  
of $L^1(G)$ is a d-ideal $\cl A.,$ and 
the map $\pi\to\pa$ is the usual one given by
$\pa(f)={\int f(x)\,\pi(x)\,d\mu(x)}.$

\thrm Theorem 4.1." If a d-ideal $\cl A.\in \cl I.(\cl R.)$ separates
$B_\sigma,$ then ${\theta:\rep\cl A.\to\cl R.}$ is surjective."
Let $\cl A.$ separate $B_\sigma.$
We first show that $\pa:\cl A.\to\cl B.(\cl H._\pi)$ is nondegenerate
for any $\pi\in\cl R..$
If $\pa$ were degenerate, there would be a nonzero $\varphi\in\cl H._\pi$ such that
$\pa(\cl A.)\varphi=0,$ i.e. 
$\omega_x\big((\psi,\,\pi(x)\varphi)\big)=0$ for all
$\psi\in\cl H._\pi,\ab\;\omega\in\cl A..$
Now $\pi\in\cl R.$ is nondegenerate, hence there is a vector $\psi\in\cl H._\pi$
such that the function 
$x\to(\psi,\,\pi(x)\varphi)$ is nonzero, and by the previous sentence 
this is in $\ker\omega$ for all $\omega\in\cl A..$
This contradicts the hypothesis that $\cl A.$ separates $B_\sigma,$ 
and thus $\pa$ is nondegenerate. We will now show that
$\pi=\theta(\pa),$ which establishes surjectivity of $\theta.$
For all $\psi,\,\varphi\in\cl H._\pi,$ $\omega\in\cl A.$ we have:
$$\eqalignno{\big(\varphi,\,\theta(\pa)(x)\,\pa(\omega)\,\psi\big)&=
\big(\varphi,\,\pa(\delta_x*\omega)\psi\big)\cr
=(\delta_x*\omega)\s y.\big((\varphi,\,\pi(y)\,\psi)\big)
&=\omega_y\big((\varphi,\,\sigma(x,y)\,\pi(xy)\,\psi)\big)\cr
=\omega_y\big((\pi(x)^*\varphi,\,\pi(y)\psi)\big)&=
\big(\varphi,\,\pi(x)\,\pa(\omega)\psi\big)\;,\cr}$$
i.e. $\theta(\pa)(x)\cdot\pa(\omega)\psi=
\pi(x)\cdot\pa(\omega)\psi$ for all $\psi\in\cl H._\pi,$ $\omega\in\cl A.\,.$
Since $\pa$ is nondegenerate, $\pa(\cl A.)\cl H._\pi$ is dense, hence
$\theta(\pa)(x)=\pi(x)$ for all $x\in G,$ which proves that
$\theta$ is surjective.

Recall that we have the canonical isometry $\iota:J_\sigma\to J_\sigma^{**}$
by $\iota(f)(\omega):=\omega(f)$ for $\omega\in J_\sigma^*,$ $f\in J_\sigma,$
and that  $J_\sigma$ is reflexive if $\iota(J_\sigma)=J_\sigma^{**}.$
If $\cl A.\subset J_\sigma^*$ is a d-ideal, we 
denote the restriction of $\iota$ by $j:J_\sigma\to\cl A.^*$ where
$j(f)(\beta):=\beta(f),$ $\beta\in\cl A.,$ $f\in J_\sigma.$ Note that 
$j$ is injective if $\cl A.$ separates $J_\sigma.$
Now even if $J_\sigma$ is not reflexive, there may still be d-ideals
$\cl A.$ such that $j(J_\sigma)=\cl A.^*,$ and we need these because:

\thrm Theorem 4.2."For a d-ideal $\cl A.\in\cl I.(\cl R.),$ the map
$\theta:\rep\cl A.\to\cl R.$ is injective with inverse map
$\pi\in\cl R.\to\pa\in\rep\cl A.$ iff  $j(J_\sigma)=\cl A.^*.$
In this case, 
$\cl A.$ is a group algebra
for ${(G,\,\theta(\rep\cl A.))}.$"
We need to prove that $\theta(\pi)\s{\cl A.}.(\omega)=\pi(\omega)$
for all $\pi\in\rep\cl A.,$ $\omega\in\cl A.$ iff
$j(J_\sigma)=\cl A.^*.$
Assume that $\theta(\pi)\s{\cl A.}.=\pi.$ Let 
$f(x):={(\varphi,\,\theta(\pi)(x)\,\psi)},$ then
$$j(f)(\omega)=\omega(f)=
\big(\varphi,\,\theta(\pi)\s{\cl A.}.(\omega)\psi\big)
=(\varphi,\pi(\omega)\psi)$$
for all $\varphi,\,\psi\in\cl H._\pi,\;\pi\in\rep\cl A.,\;\omega
\in\cl A..$ By varying the rhs over $\pi\in\rep\cl A.,$
$\varphi=\psi\in\cl H._\pi,$ we obtain all states of $\cl A.,$ and since
these span $\cl A.^*$ and $j$ is linear, it means any functional of
$\cl A.$ can be expressed as an element of $j(J_\sigma),$
i.e. $j(J_\sigma)=\cl A.^*.$\chop
Conversely, let $j(J_\sigma)=\cl A.^*.$ Now observe that
$$\eqalignno{j(f)(\omega*\beta)&=(\omega*\beta)(f)=
\omega_x\big(\beta(\lambda_xf)\big)=\omega_x\big((\delta_x*\beta)(f)\big)\cr
&=\omega_x\big(j(f)(\delta_x*\beta)\big)\qquad\quad
\forall\;\omega,\,\beta\in\cl A.,\;f\in J_\sigma\;.\cr}$$
Thus, since $j(J_\sigma)=\cl A.^*,$ we have:
$$\xi(\omega*\beta)=\omega_x\big(\xi(\delta_x*\beta)\big)
\qquad\quad\forall\;\xi\in\cl A.^*,\;\omega,\,\beta\in\cl A.\;.$$
In particular, choose $\xi(\omega)=(\varphi,\,\pi(\omega)\psi),$
$\pi\in\rep\cl A.,\;\varphi,\;\psi\in\cl H._\pi,$ then
$$\eqalignno{\big(\varphi,\,\pi(\omega*\beta)\psi\big)&=
\omega_x\big((\varphi,\,\pi(\delta_x*\beta)\psi)\big)\cr
&=\omega_x\big((\varphi,\,\theta(\pi)(x)
\,\pi(\beta)\psi\big)\cr
&=\big(\varphi,\,\theta(\pi)\s{\cl A.}.(\omega)\,\pi(\beta)\psi\big)\cr}$$
for all $\pi\in\rep\cl A.,\;\varphi,\;\psi\in\cl H._\pi,\;\omega,\;\beta
\in\cl A.\,.$ Thus 
$$\pi(\omega)\cdot\pi(\beta)\psi=\theta(\pi)\s{\cl A.}.(\omega)
\cdot\pi(\beta)\psi\;.$$
By nondegeneracy of $\pi\in\rep\cl A.$ we get 
$\pi(\omega)=\theta(\pi)\s{\cl A.}.(\omega)$ for all $\omega\in\cl A.\;.$

The condition $j(J_\sigma)=\cl A.^*$ is quite natural, if we keep in mind that
if $\cl A.$ is a group algebra, then its dual is the coefficient
space of its representation space $\cl R.,$ and the latter is 
$B_\sigma\,(=J_\sigma$ in this case by Theorem~2.1(\rn7)).

\thrm Corollary 4.3."(i) 
Any d-ideal $\cl A.\in\cl I.(\cl R.)$ which
separates $B_\sigma$ and satisfies  $j(J_\sigma)=\cl A.^*$ is a group
algebra for $(G,\,\cl R.).$\chop
(ii) Conversely 
let $\cl A.\subset J_\sigma(\cl R.)^*$ be a d-ideal 
which is a group algebra for $(G,\,\cl R.)$ where the 
$\sigma\hbox{--homomorphism}$ $\varphi:G\to UM(\cl A.)$
is obtained from the embedding of $\delta\s G.$ in the relative multiplier algebra of $\cl A..$
Then $\cl A.\in\cl I.(\cl R.),$ $\cl A.$ separates $B_\sigma$ and satisfies
$j(J_\sigma)=\cl A.^*.$"
(i) By Theorems 4.1 and 4.2, $\theta:\rep\cl A.\to\cl R.$ is bijective.\chop
(ii) If $\cl A.$ is a group algebra as stated above, then by definition
$\theta:\rep\cl A.\to\cl R.$ so $\cl A.\in\cl I.(\cl R.)\,.$
Moreover, by Proposition~1.1 each $f\in B_\sigma$ is strictly continuous,
extends uniquely by strict continuity 
to $M(\cl A.)$ and is uniquely determined by its values
on $\cl A.$ (which is strictly dense in $M(\cl A.)).$
Thus $\cl A.$ separates $B_\sigma\,.$
Finally, since $\theta$ is bijective it has inverse map
$\pi\in\cl R.\to\pi\s{\cl A.}.\in\rep\cl A.$ by:
$$\eqalignno{\big(\phi,\,\theta(\pa)(x)\,\pa(A)\psi\big)&=
\big(\phi,\,\pa(\delta_x*A)\psi\big)
=\big(\delta_x*A\big)_y\big((\phi,\,\pi(y)\psi)\big)\cr
&=\big(\delta_x)_y\big(A_z((\phi,\,\pi(y)\pi(z)\psi))\big)
=A_z\big((\phi,\,\pi(x)\pi(z)\psi)\big)\cr
&=\big((\phi,\,\pi(x)\pa(A)\psi)\big)\qquad\hbox{for all $\phi,\,\psi\in\cl H._\pi$
and $A\in\cl A..$}\cr}$$
Thus by nondegeneracy of $\pa$ it follows that $\theta(\pa)(x)=\pi(x)\,.$
Now it follows from Theorem~4.2, by the injectivity of $\theta$
that $j(J_\sigma)=\cl A.^*.$

The condition $j(J_\sigma)=\cl A.^*$ seems hard to verify in practice, so we examine 
more accessible conditions. Observe that for any $\omega,\;\beta\in J_\sigma^*$ we have
$$\leqalignno{(\omega*\beta)(f)&=\omega_x\big(\beta_y\big(\sigma(x,y)\,f(xy)\big)\big)
=\omega_x\big(\beta_y\big((\lambda_xf)(y)\big)\big)\cr
&=\omega_x\big(\big(\delta_x*\beta)(f)\big)\big)\qquad\hbox{for all $f\in J_\sigma,$}\cr
\xi(\omega*\beta)&=\omega_x\big(\xi(\delta_x*\beta)\big)\qquad
\hbox{for all $\;\xi\in \iota( J_\sigma).$}&\hbox{i.e.}\cr}$$ 
Generalising this to all $\xi\in J_\sigma^{**}$ gives a condition which is natural
for measures:
\thrm Lemma 4.4."Let $\cl A.\subset J_\sigma^*(\cl R.)$ be any d-ideal. 
If there is an $\omega\in\cl A.$ which satisfies the condition
$$\xi(\omega*\beta)=\omega_x\big(\xi(\delta_x*\beta)\big)\qquad\forall\;\xi\in \cl A.^*,\;
\beta\in\cl A.\,,\eqno{(6)}$$
then $\pi(\omega)=\theta(\pi)\s{\cl A.}.(\omega)$ for all $\pi\in\rep\cl A.\,.$
In particular, if $(6)$ holds for a dense subset of $\omega\in\cl A.,$
then ${\theta:\rep\cl A.\to\rsd}$ is injective."
Let $\xi\in \cl A.^*$ be of the form
 $\xi(\omega)=(\varphi,\,\pi(\omega)\psi),$
$\pi\in\rep\cl A.,\ab\;\varphi,\ab\;\psi\in\cl H._\pi,$ then
by Eq~(6) we find:
$$\eqalignno{\xi(\omega*\beta)=
\big(\varphi,\,\pi(\omega*\beta)\psi\big)&=
\omega_x\big((\varphi,\,\pi(\delta_x*\beta)\psi)\big)\cr
&=\omega_x\big((\varphi,\,\theta(\pi)(x)
\,\pi(\beta)\psi\big)\cr
&=\big(\varphi,\,\theta(\pi)\s{\cl A.}.(\omega)\,\pi(\beta)\psi\big)\cr}$$
for all $\pi\in\rep\cl A.,\;\varphi,\;\psi\in\cl H._\pi,\;\beta
\in\cl A.\,.$ 
By nondegeneracy of $\pi\in\rep\cl A.$ we get 
$\pi(\omega)=\theta(\pi)\s{\cl A.}.(\omega)\,.$ 
Thus if the set of $\omega\in\cl A.$ satisfying  Eq~(6) is dense in $\cl A.,$
then $\pi=\theta(\pi)\s{\cl A.}.$ for all $\pi\in\rep\cl A.\,.$

A natural class of functionals in $J_\sigma^*(\cl R.)$ to consider, are those associated with 
finite $(\sigma\hbox{-additive,}$\ab{signed}) Borel measures on $G$ according to
$\omega_\mu(f)={\int_Gf(x)\, d\mu(x),}$ with $f$ any bounded Borel function.
Note that these functionals are continous w.r.t. the supremum norm,
i.e. ${|\omega_\mu(f)|}\leq\|\omega_\mu\|\cdot\|f\|_\infty$ for $f$ bounded and Borel.
Since we need integrable maps to define such functionals, denote by $\rbg$
those representations
whose coefficient functions are Borel. So for any $\cl R.\subset\rbg,$
we can restrict the functionals $\omega_\mu$ to $B_\sigma,$ and find
$\omega_\mu\rest B_\sigma\in K_\sigma^*\subset J_\sigma^*.$
Denote the set of these functionals by $\cl M.(G)\subset K_\sigma^*(\cl R.)\,.$
Then Lemma~3.1 has a well-known extension: given $\omega_\mu$ as above,
and a $\pi\in\rbg,$ then there is
a unique operator $\pi(\omega_\mu)\in\cl B.(\cl H._\pi)$
such that $\|\pi(\omega_\mu)\|\leq\|\omega_\mu\|$ and
$$\int(\psi,\,\pi(x)\,\varphi)d\mu(x)={(\psi,\,\pi(\omega_\mu)\,\varphi)}$$
for all $\psi,\;\varphi\in\cl H._\pi.$

We will also need to integrate the map $x\to\delta_x*\beta=:h_\beta(x)\in J_\sigma^*,$
so recall the two conditions of measurability for a Banach space--valued function
w.r.t. a Borel measure $\mu,$ cf. Lemma~9, Sect~III.6.7 of Dunford and Schwartz~[DS]:\chop
(i) inverse images of Borel sets are Borel, (ii) on the complement of a null set,
the range of the function must be separable. So
 define for a given $\cl R.\subset\rbg:$
$$\eqalignno{ D_B(\cl R.):=&\Big\{\beta\in J_\sigma^*\;\Big|\;
h_\beta^{-1}(S)\;\;\hbox{is Borel when $S\subset J_\sigma^*$ is Borel,}\cr 
&\qquad\qquad\quad h_\beta(G)\;\;\hbox{is separable}\Big\}\,.\cr
F_B(\cl R.):=& D_B(\cl R.)\cap D_B(\cl R.)^*}$$
(By $D_B(\cl R.)^*$ we here mean the adjoint set in $J_\sigma^*,$ not the dual space).
If $h_\beta$ is continuous and $G$ is separable, then $\beta\in D_B(\cl R.)\,.$

\thrm Theorem 4.5."Let $\cl R.\subset\rbg,$ then\chop
(\rn1) $F_B(\cl R.)$ is a d-ideal,\chop
(\rn2) let
$\cl A.\subseteq {F_B(\cl R.)} $ be a d-ideal, and let $\omega\in
\cl M.(G)\cap\cl A..$ Then $\pi(\omega)=\theta(\pi)\s{\cl A.}.(\omega)$
for all $\pi\in\rep\cl A.,$ and hence $\theta$ is injective
on $(\rep\cl A.)\rest(\cl M.(G)\cap\cl A.).$
In particular, if $\cl A.\subset\ol{\cl M.(G)}.\cap F_B(\cl R.),$ then 
$\theta:\rep\cl A.\to\rsd$
is injective."
(\rn1) By the definition, if $\beta\in D_B(\cl R.),$ then
$h_\beta$ is measurable for any Borel measure on $J_\sigma^*.$
By Theorem~11, Sect~III.6 of Dunford and Schwartz~[DS],
such functions form a linear space and hence by Theorem~10 of the same
section in [DS], if $k:= h_\beta + h_\alpha =h_{\beta+\alpha}$ with 
$\alpha,\,\beta\in D_B(\cl R.),$ then $k^{-1}(S)$ is Borel when $S$ is Borel.
Since $k(G)\subset\ol{h_\alpha(G)+h_\beta(G)}.$ which is separable, it follows that
$k\in D_B(\cl R.),$ hence that $D_B(\cl R.)$ is a linear space.
Since convolution is continuous, the map $x\to h_\beta(x)*\alpha=h\s\beta*\alpha.(x)$
is Borel for all $\beta\in  D_B(\cl R.)$ and $\alpha\in J_\sigma^*,$
and moreover $h\s\beta*\alpha.(G)=h_\beta(G)*\alpha$ is separable.
Thus $\beta*\alpha\in D_B(\cl R.),$ i.e. $D_B(\cl R.)$ is a right ideal in 
$J_\sigma^*,$ hence an algebra. We check norm closure. Let $\{\beta_n\}\subset
D_B(\cl R.)$ be a sequence converging to $\beta\in J_\sigma^*.$
Then ${\|h_{\beta_n}(x)-h_\beta(x)\|}\to 0,$ so we obtain pointwise convergence.
Since pointwise limits of Borel maps is Borel, and 
$$h_\beta(G)\subseteq\ol{\mathop{\cup}_{n=1}^\infty h\s\beta_n.(G)}.$$
which is separable, it follows that $\beta\in  D_B(\cl R.)$
and hence that $D_B(\cl R.)$ is a Banach algebra. Since $D_B(\cl R.)$
is  a right ideal of
$J_\sigma^*,$ we have $D_B(\cl R.)*\delta_G\subseteq D_B(\cl R.).$
We also have $\delta\s G.*D_B(\cl R.)\subseteq D_B(\cl R.)$
by the following. Let $\beta\in D_B(\cl R.),$ then
$$h\s\delta_y*\beta.(x)=\delta_x*\delta_y*\beta=\sigma(x,y)\,\delta_{xy}*\beta
=\sigma(x,y)\,h_\beta(xy)$$
so by continuity of the 2-cocycle $\sigma$ and of multiplication in $G,$
it follows that this is Borel in $x.$ Moreover
$$h\s\delta_y*\beta.(G)=\set{\sigma(x,y)\,h_\beta(xy)},x\in G.\subset
\sp\big(h_\beta(G)\big),$$ which is separable.
So $\delta_y*\beta\in D_B(\cl R.),$ hence $\delta\s G.*D_B(\cl R.)\subseteq D_B(\cl R.).$
Thus $D_B(\cl R.)$ is stable under multiplication by $\delta\s G.$
and hence so is $D_B(\cl R.)^*$ which is also a Banach algebra.
Thus the Banach *-algebra $F_B(\cl R.):=D_B(\cl R.)\cap D_B(\cl R.)^*$
is also stable under multiplication by $\delta\s G.$ 
 hence is a d-ideal.\chop
(\rn2) Let $\omega\in\cl A.\cap\cl M.(G)$ with associated Borel measure $\mu.$
Now for any $\beta\in\cl A.,$ the function $x\to\delta_x*\beta\in\cl A.$
is measurable by definition of $F_B(\cl R.),$ and bounded by $\|\beta\|.$
Thus, the Bochner integral $B:={\int_G\delta_x*\beta\,d\mu(x)}$
is well--defined (cf. Chapter III~[DS]), and $B\in\cl A.\,.$
Then 
$$\xi(B)=\int_G\xi(\delta_x*\beta)\,d\mu(x)\qquad\quad\forall\;
\xi\in\cl A.^*\eqno{(7)}$$
and in particular for $\xi=j(f),\; f\in J_\sigma,$ we have
$$\eqalignno{j(f)(B)=B(f)&=\int_G(\delta_x*\beta)(f)\,d\mu(x)
=\omega_x\big((\delta_x*\beta)(f)\big)\cr
&=(\omega*\beta)(f)\qquad\quad\forall\;f\in J_\sigma\;.\cr}$$
Thus $B=\omega*\beta=\int_G\delta_x*\beta\,d\mu(x)\,,$ 
and so, using Eq.~(7) again:
$$\xi(\omega*\beta)=\int_G\xi(\delta_x*\beta)\,d\mu(x)
=\omega_x\big(\xi(\delta_x*\beta)\big)$$
for all $\xi\in\cl A.^*,\;\beta\in\cl A.,\;\omega\in
\cl M.(G)\cap\cl A..$ 
This is exactly the condition~(6) in Lemma~4.4, hence the conclusion
follows.
 
Note that we did not require that $\cl A.\in\cl I.(\cl R.),$ and so 
$\theta(\rep\cl A.)$ need not have anything to do with $\cl R..$
However, because the representations extended from $\cl A.$ to $M(\cl A.)$
are strict--strong operator continuous, it follows from the definition of
$F_B(\cl R.)$ that $\theta(\rep\cl A.)$ must consist of Borel representations.
In fact, if we only want to study convolution algebras of measures, then
it is natural to take $\cl R.=\rbg,$ and to
 analyze these algebras in $J_\sigma^*(\cl R.).$
This seems quite useful, even for locally compact groups, in that this shows we also
have a group algebra for a large set of Borel representations, which include the continuous 
representations because $L^1(G)\subset\cl M.(G)\cap F_B(\cl R.).$
Since Borel representations of Polish groups on separable Hilbert spaces must be
continuous, the discontinuous Borel representations of Polish groups must be on 
nonseparable Hilbert spaces. Such representations do occur in physics cf.~[Gr3].
\thrm Corollary 4.6." Let $\cl R.\subseteq\rbg$ and $\cl M.(G)\cap F_B(\cl R.)\not=\emptyset.$
 For any subset 
$X\subseteq\cl M.(G)\cap F_B(\cl R.)$ let 
$\cl A.(X)$ be the 
the d-ideal generated by $X$ 
and let $\cl R._X:=
\theta\big(\rep\cl A.(X)\big).$ Then \chop
(\rn1) $\cl A.(X)$ is a \sga
for the pair ${\big(G,\,\cl R._X\big)}.$\chop
(\rn2) If $\cl A.(X)$ is in $\cl I.(\cl R.)$ and
separates $B_\sigma$ then $\cl A.(X)$ is a \sga
for $(G,\,\cl R.).$"
Observe that since all $\sigma\hbox{--translations}$
of Borel measures are Borel measures, and all
$\sigma\hbox{--convolutions}$ of Borel measures
are Borel measures, ${\cl A.(X)\cap\cl  M.(G)}$
is dense in $\cl A.(X)$
(the convolutions defined here coincide with the usual ones
for measures). 
Moreover $\cl A.(X)\subset F_B(\cl R.)$
because $ F_B(\cl R.)$ is closed under algebraic operations,
multiplication by $\delta\s G.$ and w.r.t. the norm.
Thus by Theorem~4.5 the map $\theta:\rep\cl A.(X)
\to\cl R._X\subset\rsd$ is bijective, hence 
$\cl A.(X)$ is a \sga
for the pair ${\big(G,\,\cl R._X\big)}.$
By Theorem~4.1, if $\cl A.(X)\in\cl I.(\cl R.)$ separates $B_\sigma,$ then
$\cl R._X=\cl R..$ 

The subsets of $\cl M.(G)\cap F_B(\cl R.)$ behave very well, e.g. if
$X\subset Y\subseteq\cl M.(G)\cap F_B(\cl R.)$ such that for their d-ideals
$\cl A.(X)\subset\cl A.(Y),$ then by Corollary~4.6 and Proposition~1.6,
$\cl R._X$ is a direct summand of $\cl R._Y.$
In the case of $G$ locally compact and $\cl R.=\rbg,$ 
since $L^1(G)\subset \cl M.(G)\cap F_B(\cl R.),$
we know that the continuous representations must be a direct summand of
$\theta\big(\cl A.\left(\cl M.(G)\cap F_B(\cl R.)\right)\big)\,.$
\medskip

Next we would like to return to the Gelfand-Raikov problem, i.e. to consider
the existence of a \sga for $\cl R.=\rsg.$ For the rest of this section
we will maintain this choice for $\cl R.,$ unless otherwise indicated.
The first problem is to characterize $\cl I.(\cl R.)$ more explicitly.

Let us start by listing relevant structure of the embedding
$C^*(G_d)\subset M(C^*(G))$ for a locally compact group $G.$
\item{(\rn1)}
The map $G\to UM(C^*(G))$ by $x\to\delta_x$ is continuous w.r.t. the
strict topology, i.e. if $x_\nu\to x$ is a convergent net in $G,$
then ${\|(\delta_{x_\nu}-\delta_x)A\|}\rightarrow 0 \leftarrow
{\|A(\delta_{x_\nu}-\delta_x)\|}$ for all $A\in C^*(G).$
\item{(\rn2)} The action of $C^*(G_d)$ on $C^*(G)$ has cyclic elements,
in the sense that $C^*(G)=C^*\set\delta_xA,x\in G.$ for some 
$A\in C^*(G).$ 
For example we can take for $A$ any nonzero element of $C_c(G).$
\item{(\rn3)} The inverse of the extension map
$\theta:{\rm Rep}\,C^*(G)\to{\rm Rep}\,G$  
is given by $\pi\in{\rm Rep}\,G\to\pi\s{\cl A.}.\in{\rm Rep}\,C^*(G)$
where $\pi\s{\cl A.}.(f):=\int_Gf(x)\,\pi(x)\,d\mu(x),$ $f\in L^1(G)$
(cf. Lemma~3.1 for generalisation of this
representation to functionals).

\noindent Returning now to the situation where $G$ is not locally compact,
in the light of  property (i),
 we want to select
a subalgebra of $J_\sigma^*$ with good continuity properties w.r.t.
translations.
The strict continuity in property (i) is too strong
for general topological groups, so we need to consider a weaker continuity.
For any topological group $G$ we define:
$$\leqalignno{ {}\qquad\quad\qquad\cl Q._0(G)&:=\set A\in J_\sigma^*,
{\xi\big((\delta_x-\un)*\delta_y* A\big)\to 0\quad\hbox{as}\quad
x\to e\quad\forall\; y\in G,\;\xi\in J_\sigma^{**}}.&\hbox{\bf Def.}\cr
\cl Q.(G)&:=\cl Q._0(G)\cap\cl Q._0(G)^*\quad\hbox{(adjoint is meant here, not dual)}\cr
 \qquad\cl L._0(G)&:=\set A\in J_\sigma^*,
\big\|(\delta_x-\un)*\delta_y* A\big\|\to 0\quad\hbox{as}\quad
x\to e\quad\forall\; y\in G.&\hbox{\bf and}\cr
\cl L.(G)&:=\cl L._0(G)\cap\cl L._0(G)^*\quad\hbox{(adjoint here, not dual)}\cr}$$
Thus by Theorem 3.4, $A\in\cl Q._0(G)$ iff $\xi\left(A\circ
\big(\sigma(x,\, y)\,\lambda_{xy}-\lambda_y\big)\right)\to 0$
as $x\to e$ for all $y\in G$ and $\xi\in J_\sigma^{**}.$
Note that $\cl Q._0(G)\supseteq\cl L._0(G)$
and that we always have pointwise continuity
$\lim\limits_{x\to e}((\delta_x-\un)*A)(f)=0$ for all $A\in J_\sigma^*,$ 
$f\in J_\sigma$ by Theorem~2.1(\rn5).
Thus it is only possible to have $\cl Q._0(G)\not= J_\sigma^*$
if $\iota(J_\sigma)\not=J_\sigma^{**},$ i.e. if $J_\sigma$ is not reflexive.

\thrm Lemma 4.7."If  $\iota(J_\sigma)=J_\sigma^{**}$ then $G$ is discrete
and 
$J_\sigma^*$ is $\csd\,.$"
If $\iota(J_\sigma)=J_\sigma^{**}$ then by Theorem~2.1(v) we have that
$\lim\limits_{x\to e}\xi(\delta_x-\un)=0$ for all $\xi\in J_\sigma^{**}$
and hence, since for $\xi$ one can choose the coefficient functions
$\xi(A)={\big(\psi,\,\pi(A)\varphi\big)},$ $\pi\in\rep J_\sigma^*,$
it follows that $\theta(\pi)(x)$ is weak operator continuous,
hence strong operator continuous (by unitarity of $\theta(\pi)(x)).$
Thus $J_\sigma^*\in\cl I.(\cl R.)\,,$ and so we can apply Theorem~4.2
to conclude that the map $\theta:\rep J_\sigma^*\to\rsg$ is injective.
But as $\delta_G\subset J_\sigma^*,$ it follows that $J_\sigma^*$ separates
$B_\sigma,$ so by Corollary~4.3 
$J_\sigma^*$ is a group algebra
for $(G,\,\cl R.).$  This implies that
the unital subalgebra $\csd$ separates all the states of $J_\sigma^*$
hence by the Stone-Weierstrass theorem (Theorem~11.3.1 in [Di])
$\csd$ is equal to 
$J_\sigma^*.$
This has a state $\xi_0$ defined by $\xi_0(\delta_x)=1$ if $x=e,$ and 
$\xi_0(\delta_x)=0$ if $x\not=e.$ This state satisfies
the requirement that $x\to \xi_0(\delta_x)$ is continuous
iff $G$ is discrete.

Thus if $G$ is nondiscrete then  $\iota(J_\sigma)\not=J_\sigma^{**},$
and so it is possible that $\cl Q._0(G)\not= J_\sigma^*.$
We prove below that this is in fact the case.

\thrm Theorem 4.8."(\rn1) The spaces $\cl Q._0(G)$ and $\cl L._0(G)$  are norm-closed, hence so
are $\cl Q.(G)$ and $\cl L.(G).$\chop
(\rn2) If $G$ is nondiscrete, then $\delta_x\not\in\cl Q._0(G)\supseteq \cl L._0(G)$ for any $x,$
and hence $\cl Q._0(G)\not= J_\sigma^*.$\chop
(\rn3) $\cl Q._0(G)$ and  $\cl L._0(G)$ are right ideals in $J_\sigma^*,$ hence  Banach algebras.
Thus $\cl Q.(G)$ and $\cl L.(G)$  are C*-algebras.\chop
(\rn4) Both $\cl Q.(G)$ and $\cl L.(G)$  are d-ideals, i.e. $\delta\s G.$ is in their
relative multiplier algebras.\chop
(\rn5) If $G$ is locally compact, then $L^1(G)\subset\cl L.(G),$ where as usual we identify
$h\in L^1(G)$ with $\omega_h\in J_\sigma^*$ by 
$\omega_h(f):=\int h(x)\, f(x)\, d\mu(x),$ $f\in J_\sigma$
and $\mu$ the Haar measure."
 (\rn1) Consider a sequence $\{A_n\}\subset\cl Q._0(G)$ which converges in norm to
$A\in J_\sigma^*.$ Then for all $\xi\in J_\sigma^{**}:$
 $$\eqalignno{\left|\xi\big((\delta_x-\un)*\delta_y*A\big)\right|&\leq
\left|\xi\big(\delta_x*\delta_y*(A-A_n)\big)\right| +
 \left|\xi\big(\delta_x*\delta_y*A_n-\delta_y*A_n\big)\right| \cr 
 &\qquad\qquad+\left|\xi\big(\delta_y*(A_n-A)\big)\right|  \cr
 &\leq\|\xi\|\cdot\big\|\delta_x*\delta_y*(A-A_n)\big\|+
 \|\xi\|\cdot\big\|\delta_y*(A_n-A)\big\|\cr
&\qquad\qquad+\left|\xi\big((\delta_x-\un)*\delta_y*A_n \big)\right| \cr
&=2\|\xi\|\cdot\|A-A_n\|+\left|\xi\big((\delta_x-\un)*\delta_y*A_n \big)\right| \cr
&\maprightt e,x.2\|\xi\|\cdot\|A-A_n\|\maprightt \infty, n. 0\cr}$$
and thus $A\in\cl Q._0(G)$ i.e. $\cl Q._0(G)$ is norm closed.
A similar calculation establishes that $\cl L._0(G)$ is also  norm closed.\chop
(\rn2) If $\delta_x\in\cl Q._0(G)$ then by definition
$\xi\big((\delta_y-\un)*\delta_z* \delta_x\big)\to 0$
as $y\to e$ for all $z\in G$ and $\xi\in J_\sigma^{**}.$
In particular, let $z=x^{-1},$ then $\xi(\delta_y-\un)\to 0$
as $y\to e$ for all $\xi\in J_\sigma^{**}.$
However, since $\csd$ is in 
$J_\sigma^*,$
by the Hahn-Banach theorem the restriction of $J_\sigma^{**}$ to
$\csd$ is exactly the dual of $\csd,$ and we know this contains states
$\xi$ for which $\xi(\delta_y-\un)\not\to 0$  as $y\to e,$ e.g.
the state $\xi_0$ in the proof of Lemma~4.7 (since $G$ is nondiscrete).
Thus we can never have that $\delta_x\in\cl Q._0(G).$\chop
(\rn3) 
Let $A\in\cl Q._0(G)$ and $B\in J_\sigma^*,$ then for all $\xi\in  J_\sigma^{**}$
we have ${\xi\big((\delta_x-\un)*\delta_y*(A*B)\big)}=
{\xi_B\big((\delta_x-\un)*\delta_y*A\big)}$ where $\xi_B(A):={\xi(A*B)}\,.$
Obviously $\xi_B\in J_\sigma^{**}$ hence by $A\in\cl Q._0(G)$
we get that ${\xi_B\big((\delta_x-\un)*\delta_y*A\big)}\to 0$ as
$x\to e,$ and hence $A*B\in\cl Q._0(G).$ Thus $\cl Q._0(G)$ is a right ideal in $J_\sigma^*.$ 
Next let $A\in \cl L._0(G)$ and $B\in J_\sigma^*,$ then
$$\big\|(\delta_x-\un)*\delta_y*(A*B)\big\|\leq
\big\|(\delta_x-\un)*\delta_y*A\big\|\cdot\|B\|
\maprightt e,x.0$$
for all $y\in G.$ Thus $A*B\in\cl L._0(G).$
To show that $\cl L.(G)$ is a Banach *--subalgebra of $J_\sigma^*,$ note that we
already have norm--closure, and that it is closed under involution,
so it only remains to check that it is an algebra. 
Let $A,\;B\in\cl L.(G),$ hence $A,\;A^*\in
\cl L._0(G)\ni B,\;B^*.$ Since $\cl L._0(G)$ is a right ideal,
it contains $A*B,$ as well as $B^**A^*=(A*
B)^*.$ Thus $A*B\in\cl L.(G).$
By a similar argument we find that $\cl Q.(G)$ is a Banach *--algebra.\chop
(\rn4) By (\rn3) we already know that $\cl L._0(G)*\delta_x\subseteq
\cl L._0(G).$ Let $A\in\cl L._0(G),$ $z\in G,$ then
$$\eqalignno{\big\|(\delta_x-\un)*\delta_y*(\delta_z*A)\big\|&=
\big\|(\delta_x-\un)*\big(\sigma(y,z)\,\delta\s yz.\big)*A\big\|\cr
&=\big\|(\delta_x-\un)*\delta\s yz.*A\big\|\maprightt e,x.0\cr}$$
for all $y\in G.$ Thus $\delta_z*A\in\cl L._0(G),$ i.e.
$\delta_x*\cl L._0(G)\subseteq \cl L._0(G)$ for all $x\in G.$
Now let $A\in\cl L.(G)\subset\cl L._0(G),$ hence
$\delta_x*A\in\cl L._0(G),$ and
also $(\delta_x*A)^*=A^**\delta\s -x.\in\cl L._0(G)$
because $A^*\in\cl L._0(G).$ Thus $\delta_x*A\in\cl L.(G),$
and likewise $A*\delta_x\in\cl L.(G),$ hence
$\delta_x*\cl L.(G)\subseteq\cl L.(G)\supseteq\cl L.(G)*\delta_x.$
By replacing the norms $\|\cdot\|$ in the equation above by ${\big|\xi(\cdot)\big|}$
we can transcribe this argument to prove also that $\cl Q.(G)$ is a d-ideal.\chop
(\rn5) Here $G$ is locally compact. 
 Now
$$\eqalignno{\|\omega_h\|&=\sup\set\big|\omega_h(f)\big|,f\in J_\sigma,\;\|f\|_*\leq 1.\cr
&\leq\sup\set\big|\omega_h(f)\big|,f\in J_\sigma,\;\|f\|_\infty\leq 1.\cr}$$
since $\|f\|_\infty\leq\|f\|_*\,.$
Then by $\big|\omega_h(f)\big|\leq\|f\|_\infty\|h\|_1$ (where $\|\cdot\|_1$ denotes
the $L^1\hbox{--norm}),$ it follows that $\|\omega_h\|\leq\|h\|_1\,,$
and hence that the norm closure of a set in $J_\sigma^*$ contains its 
$L^1\hbox{--closure.}$
Since $\cl L.(G)$ is norm--closed,
it thus suffices to show that $C_c(G)\subset\cl L.(G).$ Let $h\in C_c(G),$
then
$$\eqalignno{\big\|(\delta_x-\un)&*\delta_y*\omega_h\big\|=
\sup\set\left|\big((\delta_x-\un)*\delta_y*\omega_h\big)(f)\right|,
\|f\|_*\leq 1.\;.&-\!(*)\cr
\hbox{Now}&\qquad\quad\Big|\big((\delta_x-\un)*\delta_y*\omega_h\big)(f)\Big|
=\big|\omega_h(\lambda_y(\lambda_x-1)f)\big|\cr
&=\Big|(\omega_h)_z\big(\big(\sigma(x,y)\,\lambda\s xy.-\lambda_y\big)
f(z)\big)\Big|\cr
&=\Big|(\omega_h)_z\big(\sigma(x,y)\,\sigma(xy,z)\,f(xyz)-\sigma(y,z)\,
f(yz)\big)\Big|\cr
&=\Big|\int_G\sigma(y,z)\Big(\sigma(x,yz)\,f(xyz)-f(yz)\Big)h(z)\,
d\mu(z)\Big|\cr
=\Big|\int_Gf(s)&\Big(\sigma(y,y^{-1}x^{-1}s)\,\sigma(x,x^{-1}s)\,
h(y^{-1}x^{-1}s)-\sigma(y,y^{-1}s)\,h(y^{-1}s)\Big)\,d\mu(s)\Big|\cr
\leq\|f\|_\infty\cdot\int_G&\Big|\ol\sigma(y^{-1},x^{-1}s).\cdot
\ol\sigma(x^{-1},s).\,h(y^{-1}x^{-1}s)-\ol\sigma(y^{-1},s).\,
h(y^{-1}s)\Big|\,d\mu(s)\cr
\leq\|f\|_*\cdot\int_G&\Big|\ol\sigma(y^{-1},x^{-1}s).\cdot
\ol\sigma(x^{-1},s).\,h(y^{-1}x^{-1}s)-\ol\sigma(y^{-1},s).\,
h(y^{-1}s)\Big|\,d\mu(s)\cr
}$$
where we made use of $\sigma(a,a^{-1}b)=\ol\sigma(a^{-1},b)..$
The integrand is bounded by $2\|h\|,$ of compact support contained in
${xy\,{\rm supp}(h)\cap y\,{\rm supp}(h)}$ and goes pointwise
to zero as $x$ approaches $e,$ independently of $f.$ Thus the Lebesgue
dominated convergence theorem applies, and we get that the last integral
goes to zero as $x$ approaches $e,$ independently of $f.$ Thus from $(*)$
we get that ${\big\|(\delta_x-\un)*\delta_y*\omega_h\big\|}\maprightt e,x.0$
for all $y\in G,$ i.e. $\omega_h\in\cl L._0(G).$ Now
$$\eqalignno{\omega_h^*(f)&=\ol\omega_h(\gamma f).
=\ol\int_Gh(x)\,{\ol f(x^{-1}).}\,d\mu(x).
=\int_G\ol h(x).\,f(x^{-1})\,d\mu(x)\cr
&=\int_G\ol h(x^{-1}).\,f(x)\,\Delta(x)\,d\mu(x)
=\omega\s\tilde{h}.(f)\cr}$$
where
$\Delta$ is the modular function of $G$
and $\widetilde{h}(x):=\ol h(x^{-1}).\Delta(x).$
As $\widetilde{h}\in C_c(G),$ it follows that $\omega^*_h\in\cl L._0(G)$
and hence $\omega_h\in\cl L.(G)$ for all $h\in C_c(G).$

To show that $\cl Q.(G)\supseteq\cl L.(G)$ is not an empty construct,
we need to give some examples of groups $G$ which are not
locally compact, with $\cl Q.(G)\not=\{0\}.$
In the next section is an example for groups
which are amenable but not locally compact
(which is a large class) and for which we show that $\cl L.(G)\not=\{0\}.$
However, an even better argument comes from the next theorem.

Recall our convention that unless otherwise specified, $G$ is
nondiscrete.

\thrm Theorem 4.9."
For a d-ideal $\cl A.$ we have that $\cl A.\in\cl I.\big(\rsg\big)$ iff
$\cl A.\subseteq\cl Q.(G)\,.$"
Let  $\cl A.\in\cl I.\big(\rsg\big)$ i.e. $\theta(\pi)\in \rsg\,.$
Now any functional $\xi\in J_\sigma^{**}$ is of the form
$\xi(B):={(\psi,\,\pi(B)\phi)},$ for $\pi\in\rep J_\sigma^*,$
$\psi,\ab\;\phi\in\cl H._\pi$ so for such a $\xi$
we have for all $A\in\cl A.:$
$$\eqalignno{\xi\big((\delta_x-\un)*\delta_y*A\big) &=
\left(\psi,\,\pi\big((\delta_x-\un)*\delta_y*A\big)\phi\right)\cr
&=\left(\psi,\,\big[\theta(\wt\pi)(x)-\un\big]\pi(\delta_y*A)\phi\right)\cr
&\maprightt e,x. 0\cr}$$  
where $\wt\pi$ is the restriction of $\pi$ to $\cl A.$ on its essential subspace,
where the latter is the closure of $\pi(\cl A.)\cl H._\pi\ni\pi(\delta_y*A)\phi\,.$
In the last step we used $\cl A.\in\cl I.\big(\rsg\big).$
Thus $\xi\big((\delta_x-\un)*\delta_y*A\big)\maprightt e,x. 0$ for all $y\in G$
and  $\xi\in J_\sigma^{**}$ i.e. $\cl A.\subseteq\cl Q.(G)\,.$\chop
Conversely, let $\cl A.\subseteq\cl Q.(G)$ and recall from the Hahn-Banach theorem
that the dual of $\cl A.$ consists of the restriction of $J_\sigma^{**}$ to
$\cl A..$ Thus ${\xi\big((\delta_x-\un)*A\big)}\to 0$ as $x\to e$ for all
$A\in\cl A.$ and $\xi\in\cl A.^*.$ By choosing coefficient functions
$\xi(A)={(\psi,\,\pi(A)\phi)},$ for $\pi\in\rep\cl A.$
we find as above that 
$$\left(\psi,\,\big[\theta(\pi)(x)-\un\big]\pi(A)\phi\right)\maprightt e,x. 0$$
for all $\psi,\;\phi$ and so $\theta(\pi)(x)$ is weak operator continuous, 
hence strong operator continuous by unitarity of $\theta(\pi)(x).$
Hence $\cl A.\in\cl I.\big(\rsg\big).$

Thus, if $\cl R.\subseteq\rsg$ and
$(G,\,\cl R.)$ has a group algebra $\cl L.,$ then $\Psi(\cl L.)\subseteq
\cl Q.(G)$ where $\Psi$ is the isometric embedding of Theorem~3.5.
So by the example below Theorem~1.5, if $\rsg$ contains an irreducible representation
then there is a nontrivial d-ideal in $\cl Q.(G)\,.$

One might be tempted to regard $\cl Q.(G)$ as a possible candidate for a group algebra,
however, in general it is too large. For example, if  $G$ is amenable, $\sigma=1,$ then all its
(two-sided) invariant means are in $\cl L.(G)\subset\cl Q.(G).$ When 
$G$ is locally compact Abelian, but not 
compact, we know it has uncountably many invariant means not identified with
elements of $L^1(G),$ cf.~[CG]. 

Now in the light of Theorem 4.9, a d-ideal $\cl A.\subset\cl Q.(G)$ will be an adequate group algebra 
for $G$ if we can show that $\theta$ is bijective. Clearly we can now use
Theorems~4.1--4.5 to find sharp conditions for this.

Note that at this point, we have obtained a distinguished subset
of representations $\theta(\rep\cl Q.(G))\subseteq\rsg$ which is the homomorphic image
of the representation theory of a C*-algebra.
By the construction preceding Prop.~1.6, we can also obtain from $\delta_G\subset M(\cl Q.(G))$
another distinguished set of representations isomorphic to
the representation theory of a C*-algebra.


\thrm Corollary 4.10." Let $\cl A.\subset\cl Q.(G)$ be a d-ideal, 
then\chop
(\rn1) 
If $\cl A.$ separates $B_\sigma$ and satisfies $j(J_\sigma)=\cl A.^*,$ then $\cl A.$
is a group algebra for $(G,\,\rsg).$\chop
(\rn2) If $G$ is separable and $\cl A.\subseteq\ol{\cl M.(G)\cap\cl L.(G)}.,$
then $\cl A.$ is a group algebra for ${(G,\,\theta(\rep\cl A.)),}$
(note that $\theta(\rep\cl A.))\subseteq\rsg).$
If in addition $\cl A.$ separates $B_\sigma,$ then $\cl A.$
is a group algebra for ${(G,\,\rsg).}$\chop
(\rn3) If $G$ is locally compact, then $\csg\cong\cl A.$
where 
 $\cl A.=\ol{\cl M.(G)\cap\cl L.(G)}..$"
 (\rn1) This follows from Corollary 4.3 and Theorem~4.9.\chop
 (\rn2) Note that $\cl R.=\rsg\subset\rbg\,,$ and that if
 $\omega\in{\cl M.(G)\cap\cl L.(G)}$ then $h_\omega(x):=\delta_x*\omega$
 is norm continuous in $x$ by definition of $\cl L.(G)\,.$
 Thus since $G$ is separable, $\omega\in F_B(\cl R.),$ 
 i.e. $\cl A.\subset  F_B(\cl R.).$ 
 Thus we can apply Theorem~4.5 and Corollary~4.6 to $\cl A..$
 By  Corollary~4.6 and Theorem~4.9 the claim above follows.\chop
(\rn3) Observe that since all $\sigma\hbox{--translations}$
of Borel measures are Borel measures, and all
$\sigma\hbox{--convolutions}$ of Borel measures
are Borel measures, ${\cl  M.(G)}\cap\cl L.(G)$
is a *-algebra stable under multiplication by $\delta_G\,,$
and so $\cl A.$ is a d-ideal.
(The convolutions defined here coincide with the usual ones
for measures). By Theorem~4.8
$L^1(G)\subseteq{\cl M.(G)\cap\cl L.(G)},$ and as it separates
$B_\sigma,$ it follows from (\rn2) that $\cl A.$ is a group
algebra for ${(G,\,\rsg).}$ By uniqueness (Theorem~1.5)
the isomorphism $\cl A.\cong\csg$ follows.

The characterization of $\csg$ in 4.10~(\rn3) above, is interesting
because it uses neither the Haar measure nor the behaviour of measures
w.r.t. compact sets, and it seems to improve the criterion in~[DvR].

Since we have an example of a group with a faithful continuous representation, but
no irreducible ones (cf.    Exmp~5.2 in [Pes]), for such a group we know
that it has no nonzero d-ideals satisfying the conditions in Corollary~4.10(i).

\beginsection 5. Example.

We want to show that $\cl L.(G)\not=\emptyset$
for some groups which are not locally compact.
Let $G$ be amenable but not locally compact
with a faithful continuous representation.
This is a large class, e.g. 
any Abelian group is amenable (and there are many examples of these
with faithful continuous representations).
For a nonabelian example, take the unitary group
 of any nuclear C*-algebra (cf. [Pa]) with the 
 relative weak topology which obviously has a faithful continuous representation
 since the C*-algebra has.

For our purposes we will take the definition of ``amenable group''
to mean that there is a left--invariant mean $n$
on $K_\sigma\supset J_\sigma.$ (This is weaker than the usual definition, since
$K_\sigma$ is uniformly continuous cf~Theorem 2.1(\rn5)).

\thrm Lemma 5.1. " Let $G$ be amenable and $\sigma=1.$ \chop
(i) There is a two-sided invariant mean $m\in K^*\subset J^*,$ i.e.
   $m=m\circ\lambda_a=m\circ\rho_a$ for all $a\in G.$\chop
(ii) Let $m\in K^*$ be a two-sided invariant mean, then $m\in\cl L.(G).$"
{\bf (i)} We adapt the usual proof, cf. 17.10 in~[HR1].\chop
Let $n\in K^*$ be a left-invariant mean and recall that
$$n^*(f):=\ol n(f^*).=\ol{n_x\big(\ol f(x^{-1}).\big)}.\quad\hbox{where}\quad
f^*(x):=\ol f(x^{-1}).\,.$$
Clearly $n^*(\un)=1$ and $n^*(f)\geq 0$ when $f\geq 0\,.$
Moreover
$$\leqalignno{\ol\big(n^*\circ\rho_a\big)(f).&=
\ol n^*_x\big(f(xa)\big).=n_x(\ol f(x^{-1}a).)=n_x\big(\ol f((a^{-1}x)^{-1}).\big)\cr
&=n_x(f^*(a^{-1}x))=n_x(f^*(x))=\ol n^*(f).\;,\cr}$$
i.e. $n^*$ is right--invariant. 
Define $m:=n*n^*,$ i.e.
$m(f):=n_x(n^*_y(f(xy)))=n^*_y(n_x(f(xy)))$ for $f\in K\,$ (cf. Theorem~3.3(i)).
Then $m\in K^*$ is positive and normalised and
$$\leqalignno{(m\circ\lambda_a)(f)&=
n^*_y(n_x((\lambda_a f)(xy)))=
n^*_y(n_x(\lambda_a(\rho_y f)(x)))\cr
&=n^*_y(n_x(\rho_y f(x)))=m(f)\;,\cr
(m\circ\rho_a)(f)&=n_x(n^*_y((\rho_a f)(xy)))
=n_x(n^*_y(\rho_a(\lambda_x f)(y)))
=n_x(n^*_y((\lambda_x f)(y)))=m(f)\,.\cr}$$
Thus $m$ is two-sided invariant.\chop
{\bf (ii)} For $m\in K^*$ a two-sided invariant mean as above,\chop
$(\delta_x-\un)*\delta_y*m(f)=m\left(\lambda\s xy.(f)-\lambda_y(f)\right)
=m(f)-m(f)=0,$ and
$$\leqalignno{(\delta_x-\un)*\delta_y*m^*(f)&=
m^*\left(\lambda\s xy.(f)-\lambda_y(f)\right) \cr
&=\ol m\big(({\lambda\s xy.}-\lambda_y)(f)^*\big).
=\ol m_z\left({\ol ({\lambda\s xy.}-\lambda_y)(f)(z^{-1}).}\right).\cr
&=\ol m_z\left({\ol f(xyz^{-1}).}-{\ol f(yz^{-1}).}\right).
=\ol m_z\left({\rho\s y^{-1}x^{-1}.f^*(z)-\rho\s y^{-1}.f^*(z)}\right).\cr
&=\ol m(f^*).-\ol m(f^*).=0\;.\cr}$$
Thus $m\in\cl L.(G)\,.$

These two--sided invariant means do not by themselves seem 
particularly useful elements of $\cl L.(G),$ 
because by the invariance they cannot separate $B_\sigma$ (cf.~Theorem~4.1),
and the only representation they can produce on $G$ via $\theta$
is the identity representation.
We define a more promising set of functionals.
\item{\bf Def.} Let $\sigma\not=1,$ and recall that $f\cdot h\in K_1$ if
$f\in K\s{\ol\sigma.}.$ and $h\in K_\sigma\,.$ Let $m\in K_1^*$
be a two--sided invariant mean of the amenable group $G,$ and define a
functional  $m^f\in K_\sigma^*$ by $m^f(h):=m(f\cdot h)=m_x(f(x)h(x)).$

\thrm Proposition 5.2." Let $G$ be amenable, and $f\in  J\s{\ol\sigma.}..$
Then $m^f\in\cl L.(G)$ for all two-sided invariant means $m\in K_1^*.$"
It suffices to prove the theorem for $f(x)=(\varphi,\,\pi(x)\psi),$
$\pi\in{\rm Rep}_{\ol\sigma.}G,$ $\varphi,\,\psi\in\cl H._\pi\,.$ Now
$$\eqalignno{&\left\|((\delta_x-\un)*\delta_y*m^f\right\|=\sup\set
\left|\left((\delta_x-\un)*\delta_y*m^f\right)(h)\right|,
h\in J_\sigma,\;\|h\|_*\leq 1.\cr
\hbox{so}\;\;& \left((\delta_x-\un)*\delta_y*m^f\right)(h)
=m\left(f\cdot\big(\sigma(x,y)\lambda\s xy.-\lambda_y\big)h\right)\cr
&=m_z\left(f(z)\cdot\big(\sigma(x,y)\,\sigma(xy,z)\,h(xyz)-\sigma(y,z)\,h(yz)\big)\right)\cr
&=m_z\left(f(y^{-1}x^{-1}z)\sigma(x,y)\,\sigma(xy,y^{-1}x^{-1}z)\,h(z)-
f(y^{-1}z)\sigma(y,y^{-1}z)\,h(z)\right)\cr
&=m_z\left(h(z)\Big(f(y^{-1}x^{-1}z)\sigma(x,y)\,\ol\sigma(y^{-1}x^{-1},z).
- f(y^{-1}z)\ol\sigma(y^{-1},z).\Big)\right)& -\!\!{\bf [1]}\cr}$$
using $\sigma(a,a^{-1}b)=\ol\sigma(a^{-1},b).\;.$ Now
$$\eqalignno{& H_{x,y}(z):= f(y^{-1}x^{-1}z)\sigma(x,y)\,\ol\sigma(y^{-1}x^{-1},z).
- f(y^{-1}z)\ol\sigma(y^{-1},z).\cr
&=\sigma(x,y)\,\ol\sigma(y^{-1}x^{-1},z).\big(\varphi,\,\pi(y^{-1}x^{-1}z)\psi\big)
-\ol\sigma(y^{-1},z).\big(\varphi,\,\pi(y^{-1}z)\psi\big)\cr
&=\sigma(x,y)\big(\pi(xy)\varphi,\,\pi(z)\psi\big) 
-\big(\pi(y)\varphi,\,\pi(z)\psi\big)&\hbox{by $\pi\in{\rm Rep}_{\ol\sigma.}G$}\cr 
&=\left(\big(\sigma(x,y)\pi(xy)-\pi(y)\big)\varphi,\,\pi(z)\psi\right)\cr
\hbox{so}\;\;&
\|H_{x,y}\|_*=\sup\set{\left|\left(\big(\sigma(x,y)\pi(xy)-\pi(y)\big)\varphi,\,\pi(A)\psi\right)\right|},
{A\in\csd,\;\|A\|\leq 1}.\,.\cr
&\leq\left\|\big(\sigma(x,y)\pi(xy)-\pi(y)\big)\varphi\right\|\cdot\|\psi\|\cr 
&\leq\big|\sigma(x,y)-1\big|\cdot\|\varphi\|\cdot\|\psi\|+
\big\|\big(\pi(xy)-\pi(y)\big)\varphi\big\|\cdot\|\psi\|\cr
&\maprightt e,x.0\qquad\forall\,y\cr
}$$ 
 Thus $\left\|H_{x,y}\right\|_*
\maprightt e,x.0$ for all $y,$ and so by 
equation~[1],\chop
$\big|\left((\delta_x-\un)*\delta_y*m^f\right)(h)\big|\leq
 \|m\|\cdot\|hH_{x,y}\|_*\leq\|m\|\cdot\|h\|_*\|H_{x,y}\|_*
\maprightt e,x.0$
for all $y.$ Thus $\big\|(\delta_x-\un)*\delta_y*m^f\big\| 
\maprightt e,x.0$ for all $y,$ i.e. $m^f\in\cl L._0(G)\,.$\chop
To prove the same for $(m^f)^*,$
$$\eqalignno{& \ol\left((\delta_x-\un)*\delta_y*(m^f)^*\right)(h).
=\ol {(m^f)^*\left(\big(\sigma(x,y)\lambda\s xy.-\lambda_y\big)h\right)}.\cr
& =m^f\left(\big(\sigma(x,y)\lambda\s xy.-\lambda_y\big)h^*\right)
 =m_z\left(f(z)\cdot\ol{\big(\sigma(x,y)\lambda\s xy.-\lambda_y\big)h(z^{-1})}.\right)\cr
&=m_z\left(f(z)\cdot\big(\ol{\sigma(x,y)\,\sigma(xy,z^{-1})\,h(xyz^{-1})}.
-\ol{\sigma(y,z^{-1})\,h(yz^{-1})}.\big)\right)\cr
&=m_z\left(f(zxy)\,\ol\sigma(x,y).\,\ol\sigma(xy,y^{-1}x^{-1}z^{-1}).
\ol h(z^{-1}).-f(zy)\ol\sigma(y,y^{-1}z^{-1}) h(z^{-1}).
\right)\cr
&=m_z\left(h^*(z)\Big(f(zxy)\ol\sigma(x,y).\,\sigma(y^{-1}x^{-1},z^{-1})
- f(zy)\sigma(y^{-1},z^{-1})\Big)\right)& -\!\!{\bf [2]}\cr}$$
by cocycle identities. Now
$$\eqalignno{& F_{x,y}(z):=(f(zxy)\ol\sigma(x,y).\,\sigma(y^{-1}x^{-1},z^{-1})
- f(zy)\sigma(y^{-1},z^{-1})\cr
&=\ol\sigma(x,y)\sigma(z,xy).\big(\varphi,\,\pi(zxy)\psi\big)-
\ol\sigma(z,y).\big(\varphi,\,\pi(zy)\psi\big)\cr
&=\ol\sigma(x,y).\,\big(\pi(z)^*\varphi,\,\pi(xy)\psi\big)
-\big(\pi(z)^*\varphi,\,\pi(y)\psi\big)\cr}$$
Thus similar to above, we find
$$\eqalignno{\|F_{x,y}\|_*
&\leq\|\varphi\|\cdot\left\|\big(\ol\sigma(x,y).\,\pi(xy)-\pi(y)\big)\psi\right\|\cr
&\leq\|\varphi\|\cdot\left\|\big(\pi(x)-\un\big)\pi(y)\psi\right\|
\maprightt e,x.0\qquad\forall\,y\cr
}$$ 
So by equation [2]:
$$\Big|\left((\delta_x-\un)*\delta_y*(m^f)^*\right)(h)\Big|
\leq\|m\|\cdot\|h\|_*\|F_{x,y}\|_*\maprightt e,x.0\qquad\forall\,y\,.$$
and thus $\Big\|(\delta_x-\un)*\delta_y*(m^f)^*\Big\|\maprightt e,x.0\qquad\forall\,y\,,$
i.e. $(m^f)^*\in\cl L._0(G)$ and hence
$m^f\in\cl L.(G)=\cl L._0(G)\cap\cl L._0(G)^*\,.$

In general there are $x\in G,$ and $f\in J\s{\ol\sigma.}.$
such that $\delta_x*m^f\not=m^f\,.$ 
Thus $G$ has sets of continuous representations which are homomorphic
images of
the representation theory of C*-algebras.

\beginsection Acknowledgement

I wish to thank Prof.~Colin Sutherland very warmly for his help, interest
and stimulating questions during
this project.

\beginsection Bibliography.

\item{[BB]} A.C. Baker, J.W. Baker: Algebras of measures on a locally compact
semigroup. J. London Math. Soc.(2) {\bf 1}, 249--259 (1969)
\item{[Ban]} W. Banaszczyk: On the existence of exotic Banach--Lie groups.
Math. Ann. {\bf 264}, 485--493 (1983)\chop
W. Banaszczyk: On the existence of commutative Banach--Lie groups which do not admit continuous
unitary representations. Colloq. Math. {\bf 52}, 113--118 (1987)\chop
W. Banaszczyk: On the existence of unitary representations of
commutative nuclear Lie groups. Studia Math. {\bf 76}, 175--181 (1983)
\item{[Bic]} K. Bichteler: A generalisation to the non-separable case
of Takesaki's duality theorem for C*-algebras. Invent. Math. {\bf 9}, 
89--98 (1969)
\item{[Bus]} R.C. Busby: Double centralizers and extensions of C*-algebras.
Trans. Amer. math. Soc. {\bf 132}, 79--99 (1968)
\item{[CG]} C. Chou: On topologically invariant means on a locally compact group.
Trans. Amer. Math. Soc. {\bf 151}, 443--456 (1970).\chop
 E. Granirer: On the
invariant mean on topological semigroups and on topological groups. Pac. J. Math.
{\bf 15}, 107 (1965).
\item{[Di]} J. Dixmier; C*--algebras, North--Holland, Amsterdam, 1977.
\item{[Do]} R.G. Douglas: Generalized group algebras. Illinois J. Math.
{\bf 10}, 309--321 (1966)
\item{[DS]} N. Dunford, J. Schwartz, Linear operators, Vol I. (Classics edition)
Wiley--Interscience, New York, London, Sydney, Toronto (1988) 
\item{[DvR]} H.A.M. Dzinotyiweyi, A.C.M. van Rooij:  A characterization of the group algebra of 
noncompact locally compact topological group. Quart. J. Math. Oxford Ser. (2)
{\bf 41} no 161, 15--20 (1990)
\item{[Dz]} H.A.M. Dzinotyiweyi: The analogue of the group algebra for topological
semigroups. Pitman Advanced Publishing Program, Boston, London, Melbourne (1984)
\item{[Gre]} F.P. Greenleaf: Characterization of group algebras in terms of
their translation operators. Pacific J. Math. {\bf 18}, 243--276 (1966)
\item{[GR]} I.M. Gelfand, D.A. Raikov: Irreducible representations
of locally bicompact groups, Mat. Sbornik N.S. {\bf 13}, 301--316 (1943)
(Russian). English translation: Transl. Amer. Math. Soc. (II Ser.)
{\bf 36} 1--15 (1964)
\item{[Gr1]} H. Grundling, A group algebra for inductive limit groups.
Continuity Problems of the Canonical Commutation Relations.
Acta Applicandae Math. {\bf 46}, 107--145, 1997. 
\item{[Gr2]} H. Grundling, Host Algebras: Generalising Group Algebras.
Submitted. Preprint at http://arXiv.org/abs/math/?0007112.
\item{[Gr3]}  H. Grundling, C.A. Hurst: A note on regular states and
supplementary conditions. Lett. Math. Phys. {\bf 15}, 205--212 (1988)
\item{[HR1]} E. Hewitt, K.A. Ross: Abstract Harmonic Analysis I.
Berlin, Goettingen, Heidelberg: Springer--Verlag 1963
\item{[Ma]} G.W. Mackey: The theory of unitary group representations.
Chicago: The University of Chicago Press 1976
\item{[Man]}
Manuceau, J.: C*--alg\`ebre de relations de commutation. Ann.
  Inst. H. Poincar\'e {\bf 8}, 139--161 (1968)
\item{[Pa]} A. Paterson: Nuclear C*-algebras have amenable
unitary groups. Proc. Am. Math. Soc. {\bf 114}, 719-721 (1992)
\item{[Pes]}
V. Pestov: Abelian topological groups without irreducible Banach
representations. Abelian groups, module theory and topology
(Padua, 1997), 343--349, Lecture Notes in Pure and Appl. Math. 201,
Dekker, New York 1998.
\item{[PR]} J.A. Packer, I. Raeburn: Twisted crossed products of C*--algebras.
Math. Proc. Camb. Phil. Soc. {\bf 106}, 293--311 (1989)
\item{[Py]} J.S. Pym: The convolution of linear functionals.
  Proc. Lond. Math. Soc. {\bf 14}, 431 (1964)
\item{[Tak]} M. Takesaki: A duality in the representation theory of 
C*-algebras. Ann. Math. {\bf 85}, 37---382 (1967)
\item{[Wa]} M.E. Walter: W*--algebras and nonabelian harmonic analysis.
J. Funct. Anal. {\bf 11,} 17--38 (1972)
\item{[Wh]} R.F. Wheeler: A survey of Baire measures and strict topologies.
Expositiones Math. {\bf 2}, 97--190 (1983)
\item{[WO]} N.E. Wegge-Olsen: K-theory and C*-algebras. Oxford University Press,
Oxford, New York, Tokyo (1993)
\item{[Wor]} S.L. Woronowicz: C*-algebras generated by unbounded elements.
Rev. Math. Phys. {\bf 7}, 481--521 (1995)

\bye